\documentclass[letterpaper,oneside]{amsart} 

\usepackage[T1]{fontenc}
\usepackage{color,soul}
\usepackage{amsrefs}
\usepackage[pagewise]{lineno}
\usepackage{xcolor}
\usepackage[utf8]{inputenc} 
\usepackage{comment}
\usepackage{mathtools} 
\usepackage{colonequals}
\usepackage{amsmath}
\usepackage{amssymb}
\usepackage{tikz}
\usepackage{mathrsfs}
\usepackage{enumitem}
\usepackage[foot]{amsaddr}
\usepackage{mathtools}

\usepackage{mleftright}

\usepackage[english]{babel} 
\usepackage{tikz}
\usepackage{pdfpages}
\usepackage{thmtools} 
\usepackage{amssymb}
\usepackage{hyperref}
\usepackage[all]{xy} 
\usepackage[normalem]{ulem}
\usepackage{todonotes}
\setuptodonotes{inline}

\newcommand{\Id}{\operatorname{Id}}

\newcommand{\Absurd}{(\rightarrow\leftarrow)}

\title{Universal coverings for limit and pseudotopological spaces}
\author{Jonathan Treviño-Marroquín}
\thanks{This material is based in part upon work supported by the US National Science Foundation
under Grant No. DMS-1928930 while the author participated in a program supported by the
Mathematical Sciences Research Institute. The program was held in the summer of 2022 in partnership with the the Universidad Nacional Aut\'onoma de M\'exico. This work was also supported
by the grant N62909-19-1-2134 from the US Office of Naval Research Global and the Southern Office of Aerospace Research and Development
of the US Air Force Office of Scientific Research. The author was also supported by the
CONACYT postgraduate studies scholarship number 839062.}

\theoremstyle{plain}
\newtheorem{teorema}{Theorem}[subsection]
\newtheorem{proposicion}[teorema]{Proposition}
\newtheorem{corolario}[teorema]{Corollary}
\newtheorem{lema}[teorema]{Lemma}

\theoremstyle{definition}
\newtheorem{definicion}[teorema]{Definition}

\newtheorem{ejemplo}{Example}[teorema]

\declaretheoremstyle[
qed=\qedsymbol
]{mystyle}

\begin{document}

\begin{abstract}
Limit and Pseudotopological spaces are two generalizations of topological spaces which are defined by indicating what filters converge under some axioms. In this article, we introduce covering spaces and set forth some necessary conditions for a construction for a universal covering space.
\end{abstract}

\maketitle

\tableofcontents 



\section{Introduction}
\pagenumbering{arabic}\label{cap.introduccion}
The study of homotopical invariants of discrete spaces seen through a variety of categories has become a topic of significant recent interest, due to a number of applications in different areas. One recent approach has been to study discrete spaces as Pseudotopological spaces (or Choquet spaces) \cite{Beattie_Butzmann_2002, Preuss_2002, Rieser_arXiv_2022} which is the cartesian closed hull of pretopological spaces (or \v{C}ech closure spaces) \cite{Beattie_Butzmann_2002}. The homotopy theory of pseudotopological spaces has been proposed as a possibility for the homotopy theory of points clouds and graphs, with the advantage that, in pseudotopological spaces, one has non-necessarily-trivial maps to and from topological spaces, which aids in the comparison of the invariants of discrete structures and topological spaces. A property that is significant in pseudotopological spaces is the existence of all the small limits and colimits, which is an important advantage over other categories because it opens the significantly more structure. Indeed, it has been recently shown that there is a natural model category structure in pseudotopological spaces extending the Quillen model structure on topological spaces \cite{Rieser_arXiv_2022}. Related works in other categories are $A$-homotopy \cite{Babson_etal_2006}, $\times$-homotopy \cite{Dochtermann_2009, Chih_Scull_2021} and the developing of different homotopies in \v{C}ech closure spaces \cite{Bubenik_Milicevic_2021b,Rieser_2021}, as well as homotopy in digraphs \cite{Grigoryan_etal_2014, kapulkin2023fundamental}. Most of these tools look for homotopical invariants which allow us to study graphs and simplicial complexes for applications to combinatorics \cite{BarceloSeversWhite_2011_kParabolic, BarceloSeversWhite_2013_TheDiscFund}, or in the development of topological data analysis (to see the relation with this topic see for example \cite{Carlsson_2009, Carlsson_et_al_2012, Demaria_1987,Rieser_2021}). 
 
In this work, we develop covering space theory and construct universal covers in this context, although actually work in limit spaces, since no additional problems are encountered there. We work with the classical definition of homotopy (the interval is the topological space $[0,1]$, denoted by $I$, and the product is the categorical product, denoted by $\times$). In related work, universal covers have been realized in $A$-theory for  graphs in \cite{kapulkin2023fundamental}.
The main complication in the construction arises from the loss of the neigborhood filter when passing from topological or pretopological spaces to pseudotopological spaces. These filters are essential for the definition of  open (or interior) covers. Nonetheless, by using local covering systems (\ref{def:local covering system}) in place of open (or interior) coverings, we may successfully handle the lack of (good) open sets in pseudotopological spaces. This change also required a generalization of the notions locally and semi-locally connectedness. Once all of this is established, we may complete the essential step of lifting local covering systems from $X$ to one in the path space $P(X,x_0)$ (\autoref{def:PathSpace}).

In order to define the local path connected property, we define a base for local coverings (\autoref{def:LocalCoveringSystemBase}). We then show that the necessary conditions to obtain a universal covering are correspond to those in the topological case and a universal covering must be simply connected. We use this to construct many examples of simply connected limit spaces. As mentioned above, the article is written for limit spaces, which contains pseudotopological spaces as a full subcategory (even a reflective one). However, in the last section we observe that the proof in the pseudotopological spaces are analogous, and their main difference is in what category we take the quotient. In addition, we conclude  proving that if $X$ is a pseudotopological space, then the resulting universal covering
	spaces are exactly the same when constructed in either limit or pseudotopological spaces.

\section{Limit Spaces}
\label{sec:Limit Spaces}

Limits spaces and its continuous maps form a category which contains topological spaces, graphs and metric spaces with a privileged scale. It is also cartesian closed, which makes the evaluation map continuous.

In this section, we introduce limit spaces and we give the necessary background for the sections about covering spaces. Note that limit spaces are some times called convergence spaces, as in \cite{Beattie_Butzmann_2002}. We follow the terminology in \cite{Preuss_2002}.

\subsection{Basic notions and examples}
\label{sub:Basic notions and examples}

Recall that a \emph{filter} on a set $X$ is a non-empty collection of subsets of $X$ which does not contain the empty set and is closed under supersets and finite intersections. A subset $\mathcal{B}$ of a filter $\mathcal{F}$ is called a (filter) \emph{basis} of $\mathcal{F}$ and $\mathcal{F}$ a \emph{filter generated by} $\mathcal{B}$ if and only if each set in $\mathcal{F}$ contains a set of $\mathcal{B}$. We then write $\mathcal{F}=[\mathcal{B}]$ or $\mathcal{F}=[\mathcal{B}]_X$. Also, we write $[A]$ for $[\{A\}]$ and $[x]$ for $[\{x\}]$ where $A\subset X$ and $x\in X$, respectively.

If $\mathcal{F}$ and $\mathcal{G}$ are filters on $X$, then $\mathcal{F}$ is said to be \emph{finer} than $\mathcal{G}$ and $\mathcal{G}$ \emph{coarser} than $\mathcal{F}$ when $\mathcal{G}\subset\mathcal{F}$. Let $\mathcal{F}$ be a filter which is not properly contained in any other filter. We say that $\mathcal{F}$ is an ultrafilter. Remember that for every filter $F$ exists an ultrafilter $\mathcal{G}$ such that $\mathcal{G}\supset \mathcal{F}$ via the axiom of choice (\cite{Jech_2003}, 7.5, Tarski).

Let $f:X\rightarrow Y$ be a mapping and $\mathcal{F}$ be a filter on $X$. We first denote $f(\mathcal{F})$ the filter generated by the $\{f(A):A\in\mathcal{F}\}$.

\begin{definicion}[Limit Space; \cite{Beattie_Butzmann_2002}, 1.1.1]
\label{def:Limit Space}
Let $X$ be a set. A mapping $\lambda$ from $X$ into the power set of the set of all filters on $X$ is called a \emph{limit structure} on $X$ and $(X,\lambda)$ a \emph{limit space} if and only if the following hold for all $x\in X$:
\begin{enumerate}[label=(\roman*)]
	\item $[x]$ belongs to $\lambda(x)$.
	\item For all filter $\mathcal{F},\mathcal{G}\in\lambda(x)$ the intersection $\mathcal{F}\cap\mathcal{G}\coloneqq\{U\mid U\in\mathcal{F},U\in\mathcal{G}\}$ belongs to $\lambda(x)$.
	\item If $\mathcal{F}\in\lambda(x)$ then $\mathcal{G}\in\lambda(x)$ for all filters $\mathcal{G}$ on $X$ which are finer than $\mathcal{F}$.
\end{enumerate}
Whenever it is clear from context, we just write $X$ without the limit structure. We also write $\mathcal{F}\rightarrow x$ for $\mathcal{F}\in \lambda(x)$.
\end{definicion}

We now define morphisms for these objects. Observe that the idea comes directly from topological continuity.

\begin{definicion}[\cite{Beattie_Butzmann_2002}, 1.1.3]
Let $X$ and $Y$ be limit spaces. A mapping $f:X\rightarrow Y$ is called \emph{continuous} at a point $x\in X$ if $f(\mathcal{F})\rightarrow f(x)$ in $Y$ whenever $\mathcal{F}\rightarrow x$ in $X$. The mapping $f$ is called \emph{continuous} if it is continuous at every point of $X$, and $f$ is a \emph{homemomorphism} if it is bijective and both $f$ and $f^{-1}$ are continuous.
\end{definicion}

Limit spaces has its own definition of neighborhood. We define a \emph{neighborhood system} for a limit space as the following filter $\mathcal{U}(x):=\cap\{\mathcal{F}:\mathcal{F}\rightarrow x\}$ (\cite{Beattie_Butzmann_2002}, definition 1.3.1).

\begin{ejemplo}
\label{exam:Limit_spaces}.
\begin{enumerate}[label=\roman*)]
\item Since a topology determines that every filter finer than the neighborhood filter at $x$ converges, every topological space has a natural limit space induced by the convergence of its filters. In particular, we denote by $I:=[0,1]$ with the limit structure given by its metric topology.
\item We emphasize that two different limit spaces can have the same neighborhood system for a point $x$; as a simple example take a non-finite set $X$ and the following limit spaces: one where $\{X\}\rightarrow x$, other one where $[\{x_1,\ldots,x_n\}]_X\rightarrow x$ for every $n\in\mathbb{N}$. $\mathcal{U}(x)=\{X\}$.
\item \label{def:pseudotopological_Spaces} (\cite{Beattie_Butzmann_2002}, 1.3.23) A limit space $X$ is called a \emph{pseudotopological} (or \emph{pseudotopological}) space if $\mathcal{F}\rightarrow x$ in $X$ whenever every ultrafilrer $\mathcal{G}$ finer than $\mathcal{F}$ converges to $X$.
\item \label{def:Pretopological_space} (\cite{Beattie_Butzmann_2002}, 1.3.15)
A limit space $X$ is called \emph{pretopological} if $\mathcal{U}(x)$ converges to $x$ in $X$ for every $x\in X$, i.e., if the neighborhood filter of each point converges to this point.
\end{enumerate}
\end{ejemplo}

\ref{def:pseudotopological_Spaces} and \ref{def:Pretopological_space}, of \ref{exam:Limit_spaces}, and continuous maps are both full subcategories of limit spaces. Observe that \ref{def:Pretopological_space} spaces are \v{C}ech closure spaces \cite{Cech_1966,Rieser_2021}.

\begin{ejemplo}
\label{exam:ChoquetNonCech}
Pseudotopological and pretopological spaces are clearly not the same. As simple example take a fixed $r\leq 0$ and $X=S^{1}$ such that a filter $\mathcal{F}$ converges to $x$ if an only if exists $n\in\mathbb{N}$ such that $\mathcal{F}=[\{x_1,\ldots,x_n\}]_X$ and $d(x,x_i)\leq r$ for every $i\in \{1,\ldots,n\}$.

This space is a limit space. $[x]\rightarrow x$ for every $x\in X$ by definition. If $\mathcal{F},\mathcal{G}\in\lambda(x)$, then $\mathcal{F}=[\{x_1,\dots,x_n\}]$, $\mathcal{G}=[\{x'_1,\ldots,x'_m\}]$, $d(x_i,x)\leq r$ for every $i\in\{1,\ldots,n\}$ and $d(x'_i,x)\leq r$ for every $i\in\{1,\ldots,m\}$; thus $\mathcal{F}\cap\mathcal{G}=[\{x_1,\ldots,x_n,x'_1,\ldots,x'_m\}]\in\lambda(x)$. Finally, if $\mathcal{F}\in\lambda(x)$, then $\mathcal{F}=[\{x_1,\ldots,x_n\}]$ with $d(x_i,x)\leq r$; thus if $\mathcal{F}'$ is finer than $\mathcal{F}$, $\mathcal{F}'=[A]$ with $A\subset \{x_1,\ldots,x_n\}$.

This space is pseudotopological. If $\mathcal{F}$ is an ultrafilter contained in $[\{x_1,\ldots,x_n\}]$, then $\mathcal{F}=[x_i]$ for some $i\in\{1,\ldots,n\}$.

This space is not pretopological. Let $x\in X$, then $A\in \mathcal{U}(x)=\bigcap\limits_{\mathcal{F}\rightarrow x}\mathcal{F}$ has to contain every $y\in X$ such that $d(y,x)\leq r$. Thus, $A\supset \{y\in S_1: d(x,y)\leq r\}$, i.e. $\mathcal{U}(x)\supset [\{y\in S_1: d(x,y)\leq r\}]$ and their elements have an infinite amount of points.
\end{ejemplo}

We have the following simple result which is very useful. With this lemma, we are able to replace any filter $\mathcal{F}\rightarrow x$ with a filter $\mathcal{F}'\rightarrow x$ such that every set in $\mathcal{F}'$ contains $x$.

Note that the limit spaces in \ref{exam:Limit_spaces}
have many filters $\mathcal{F}\rightarrow x$ with sets $U\in\mathcal{F}$ that do not contain $x$. To mention one, $\mathcal{F}=[\{ [0,1]\cap (\frac{1}{2}-\frac{1}{n},\frac{1}{2}): n\in\mathbb{N},n\geq 2\}]\rightarrow \frac{1}{2}$ and the elements $(\frac{1}{2}-\frac{1}{n},\frac{1}{2})$ does not include $\frac{1}{2}$.

\begin{lema}
\label{lem:Filter with point}
Let $X$ be a limit space. $\mathcal{F}\rightarrow x$ if and only if $\mathcal{F}\cap [x] \rightarrow x$.
\end{lema}

\begin{proof}
Let $X$  be a limit space and $\mathcal{F}$ a filter on $X$.

$(\Rightarrow)$ Assume that $\mathcal{F}\rightarrow x$. Since $[x]\rightarrow x$ and filters are closed under intersection in limit spaces, then we have that $\mathcal{F}\cap [x]\rightarrow x$.

$(\Leftarrow)$ Now assume that $\mathcal{F}\cap [x]\rightarrow x$. Observe that $\mathcal{F}\supset \mathcal{F}\cap [x]$. Thus $\mathcal{F}\rightarrow x$.
\end{proof}

Since \ref{exam:ChoquetNonCech}, we observe that using just neighborhood systems does not describe all of the information in limit spaces, i.e. we lost some convergent filters. Thus, taking just open covers is not enough for limit spaces. In order to construct covering spaces, we have the following definition which generalizes an open cover.

\begin{definicion}[Local Covering System; \cite{Beattie_Butzmann_2002}, 1.3.28]
\label{def:local covering system}
Let $X$ be a limit space and $x\in X$. A \emph{local covering system} (or \emph{local covering} or \emph{local cover}) of $X$ at $x$ is a collection $\mathcal{C}$ of subsets of $X$ such that $\mathcal{C}\cap \mathcal{F}\neq\emptyset$ for each $\mathcal{F}\rightarrow x$ in $X$.

Let $A$ be a subset of $X$. A \emph{covering system} of $A$ is meant a collection $\mathcal{C}$ of subsets of $X$ which is a local covering system of $X$ at each point of $A$. Whenever we omit the set, we will assume that $A=X$.
\end{definicion}

\subsection{New spaces generated by old ones}
\label{sub:New spaces}

Topological constructs (see \cite{Preuss_2002}) admit initial and final limit structures, which enables us to build new space from old ones; for example product structures, subspaces, and quotients. We give the definition of them in this section, as well as the definition of the continuous limit structure or spaces of continuous functions $\mathcal{C}(X,Y)$.

\begin{definicion}[Topological constructs; \cite{Preuss_2002}, 1.1.1 and 1.1.2]
\label{def:topological_construct}
By a \emph{construct} we mean a category $\mathcal{C}$ whose objects are structured sets, i.e. pairs $(X,\xi)$ where $X$ is a set and $\xi$ a $\mathcal{C}$-structure on $X$, whose morphism $f:(X,\xi)\rightarrow (Y,\eta)$ are suitable maps between $X$ and $Y$ and whose composition law is the usual composition map.

A construct $\mathcal{C}$ is called \emph{topological} if and only if it satisfies the following conditions:
\begin{enumerate}[label = (\arabic*)]
\item \emph{Existence of initial structures}:
For any set $X$, any family $((X_i,\xi_i))_{i\in I}$ of $\mathcal{C}.$-objects indexed by a class $I$ and any family $(f_i:X\rightarrow X_i)_{i\in I}$ of maps indexed by $I$ there exists a unique $\mathcal{C}$-structure $\xi$ on $X$ which is \emph{initial} with respect to $(X,f_i,(X_i,\xi_i),I)$, i.e. such that for any $\mathcal{C}$-object $(Y,\eta)$ a map $g:(Y,\eta)\rightarrow (X,\xi)$ is a $\mathcal{C}-morphism$ if and only if for every $i\in I$ the composite map $f_i\circ g:(Y,\eta)\rightarrow (X_i,\xi_i)$ is a $\mathcal{C}$-morphism.

\item For any set $X$, the class $\{(Y,\eta)\in Ob(\mathcal{C}): X=Y\}$ of all $\mathcal{C}$-objects with underlying set $X$ is set.

\item For any set $X$ with cardinality at mos one, there exists exactly one $\mathcal{C}$ object with underlying set $X$ (i.e. there exists exactly one $\mathcal{C}$-structure on $X$).
\end{enumerate}
\end{definicion}

\begin{teorema}[Existence of final structures\cite{Preuss_2002}, 1.2.1.1]
Let $\mathcal{C}$ a construct. Then the following are equivalent:
\begin{enumerate}[label = (\alph*)]
\item Existence of initial structures.
\item \emph{Existence of final structures}:
For any set $X$, any family $((X_i,\xi_i))_{i\in I}$ of $\mathcal{C}$-objects indexed by some class $I$ and any family $(f_i:X_i\rightarrow X)_{i\in I}$ of maps indexed by $I$ there exists a unique $\mathcal{C}$-structure $\xi$ on $X$ which is \emph{final} with respect to $((X_i,\xi_i),f_i,X,I)$, i.e. such that for any $\mathcal{C}$-object $(Y,\eta)$ a map $g:(X,\xi)\rightarrow (Y,\eta)$ is a $\mathcal{C}$-morphism if and only if for every $i\in I$ the composite map $g\circ f_i: (X_i,\xi_i)\rightarrow (Y,\eta)$ is a $\mathcal{C}$-morphism.
\end{enumerate}
\end{teorema}

\begin{teorema}[\cite{Preuss_2002}, 2.2.12 and 2.3.1.5]
Limit spaces, pseudotopological spaces and Pretopological spaces are topological constructs.
\end{teorema}

Observe that in the following definitions we mention \textbf{limit} structure. The constructions for other topological constructs might be different, even when they satisfy the same property. For more information about these differences see \cite{Preuss_2002}.

\begin{definicion}[Initial limit structure; \cite{Beattie_Butzmann_2002}, 1.2.1]
\label{def:Initial limit structure}
Let $X$ be a set, $(X_i)_{i\in I}$ a collection of convergence spaces and, for each $i\in I$, $f_i:X\rightarrow X_i$ a mapping. A filter $\mathcal{F}$ converges to $x$ in the \emph{initial limit structure} on $X$ with respect to $(f_i)_{i\in I}$ if and only if for each $i\in I$, $f_i(\mathcal{F})\rightarrow f_i(x)$ in $X_i$.
\end{definicion}

Two constructions which follow from from the initial structure are product and subspace structures.

\begin{definicion}[Product limit structure; \cite{Beattie_Butzmann_2002}, example 1.2.2]
\label{def:Product limit structure}
Let $(X_i)_{i\in I}$ be a collection of limit spaces and let $\prod_{i\in I}X_i$ be the product set of the $X_i$. The \emph{product limit structure} on $\prod_{i\in I}X_i$ is the initial limit structure with respect to the projection mappings $(p_i:\prod_{j\in I}X_j\rightarrow X_i)_{i\in I}$. The resulting limit space is said to be the \emph{product} of the $(X_i)_{i\in I}$.
\end{definicion}

\begin{definicion}[Subspace limit structure; \cite{Beattie_Butzmann_2002}, example 1.2.2]
\label{def:Subspace limit structure}
Let $X$ be a limit space and let $M$ be a subset of $X$. The \emph{subspace limit structure} on $M$ is the initial limit structure with respect to the inclusion mapping $e:M\rightarrow X$.
\end{definicion}

In the following definitions, we are using the existence of final limit structures. Nonetheless, we give the following characterization because is not our intention to discuss final structures neither in quotient in this part of the article \cite{Preuss_2002}.

\begin{definicion}[Quotient limit structure; \cite{Beattie_Butzmann_2002}, example 1.2.10]
\label{def:Quotients}
Let $(X,\lambda)$ be a limit space, $Y$ a set and $q:X\rightarrow Y$ a surjection. By the \emph{quotient limit structure} on $Y$ we meant that $\mathcal{F}\in\lambda_q(y)$ if and only if there are points $x_1,\ldots,x_n$ in $X$ and for each $k\in \{1,\ldots,n\}$ a filter $\mathcal{F}_k\in\lambda(x)$ such that $q(x_k)=y$ for all $k$ and $q(\mathcal{F}_1)\cap\ldots\cap q(\mathcal{F}_n)\subset \mathcal{F}$.
\end{definicion}

Given limit spaces $X$ and $Y,$ we denote the set of all continuous mappings from $X$ to $Y$ by $\mathcal{C}(X,Y)$ and we denote by $\omega_{X,Y}:\mathcal{C}(X,Y)\times X\rightarrow Y$ the \emph{evaluation mapping}, i.e. $\omega_{X,Y}(f,x)=f(x)$ for all $f\in\mathcal{C}(X,Y)$ and all $x\in X$.

\begin{definicion}[Continuous limit structure; \cite{Beattie_Butzmann_2002}, 1.1.5]
\label{def:ContinuousLimitStructure}
Let $X$ and $Y$ be limit spaces. Then the \emph{continuous limit structure} $\lambda_c$ on $\mathcal{C}(X,Y)$ is defined as follows: $\mathcal{H}\rightarrow f\in (\mathcal{C}(X,Y),\lambda_c)$ if and only if $\omega_{X,Y}(\mathcal{H}\times \mathcal{F})\rightarrow f(x)$ for all $x\in X$ and all $\mathcal{F}\rightarrow x$ in $X$, where $\mathcal{H}\times \mathcal{F}\coloneqq \{ A\times B: A\in\mathcal{H},B\in\mathcal{F}\}$.
\end{definicion}

In \autoref{sec:Covering Spaces}, we are especially interested in based spaces and paths, so we will make the following definition.

\begin{definicion}
\label{def:PathSpace}
Let $Y$ be a limit space $y_0\in Y$. We call $P(Y,y_0):=\{\gamma\in C(I,Y):\gamma(0)=y_0\}$ the \emph{path space} and we assign it the subspace limit structure of $\mathcal{C}(I,Y)$.
\end{definicion}

\subsection{Homotopy}
\label{sub:Homotopy}

Here we define what homotopy for limit spaces is in this paper. 

\begin{definicion}
\label{def:homotopy functions}
Let $X$ and $Y$ limit spaces and let $g,f:X\rightarrow Y$ be mappings. We say that $f$ and $g$ are homotopic, written $f\simeq g$, if there exists a continuous $H:X\times I\rightarrow Y$ such that $H(x,0)=f(x)$ and $H(x,1)=g(x)$.
\end{definicion}

We know that homotopy is an equivalence relation. We only need to define the first homotopy group for this paper. See \cite{Rieser_arXiv_2022} for the definition of higher homotopy groups.

\begin{definicion}
\label{def:First homotopy group}
Let $X$ be a limit space and $x\in X$. We define $\pi(X,x):=\{[g]:g:(S^1,*)\rightarrow (X,x)\}$ called the \emph{first homotopy group}.
\end{definicion}

\section{Covering Spaces and Fibrations}
\label{sec:Covering Spaces}

In this section we talk about connectedness, covering spaces and fibrations. We study covering spaces as \cite{Spanier_1966}, by restructuring some concepts and proofs to consider filters other than neighborhood filter.

We introduce our proposal. To achieve our objective, we modify the concept of local connectedness in such a way that it does not depend on an open neighborhood basis. We observe in \ref{def:LocalCoveringSystemBase} that is not enough to take a basis of every convergent filter, but a basis of coarser convergent filters which allows many connected subsets as generators.

The results of our study are a possible candidate for the universal covering space and the lifting theorem for limit spaces.

\subsection{Covering Spaces}
\label{sub:CoveringSpaces}

We start with the definition of connectedness.

\begin{definicion}
\label{def:connected}
Let $X$ be a limit space. We say that $X$ is \emph{disconnected} if there exist non-empty subsets $A$ and $B$ of $X$ such that:
\begin{itemize}
	\item $X=A\cup B$ and $A\cap B=\emptyset$.
	\item there exists a covering system $\mathcal{U}$ such that for every $U\in\mathcal{U}$ either $U\subset A$ or $U\subset B$.
\end{itemize}
Otherwise, $X$ is said to be \emph{connected}.

We call $C$ a (connected) component if $C$ is connected with the subspace limit structure and $C=C'$ if $C'\supset C$ is connected.
\end{definicion}

Observe that disconnectedness condition is equivalent to saying that exists $A$ and $B$ subsets of $X$ such that $X=A\cup B$, $A\cap B=\varnothing$ and for every $\mathcal{F}\rightarrow x$ either $A\in\mathcal{F}$ or $B\in\mathcal{F}$.

\begin{ejemplo}
\hspace{1mm}
\begin{enumerate}
\item Every connected topological space is connected.
\item \ref{exam:ChoquetNonCech} is connected.

Assume that that space is disconnected. Then there exist $A,B\subset X$ non-empty such that $A\cap B=\varnothing$ and $A\cup B=X$ and $\{A,B\}$ is a local covering system.

Suppose that there exists $a\in A$ such that there exists $b\in B$ with $d(a,b)\leq r$, then $[\{a,b\}]_X\rightarrow a,b$ $(\rightarrow \leftarrow)$. This contradicts that $\{A,B\}$ is a local covering system.

Now suppose that for every $a\in A$ there is not exist $b\in B$ with $d(a,b)\leq r$. This implies that $B=\varnothing$ $\Absurd$. $B$ has to be non-empty. Thus, $X$ is connected.
\end{enumerate}
\end{ejemplo}

The following is a useful equivalence. Under some assumptions, we can use it to construct a path between two points in a connected space.

\begin{proposicion}
\label{lem:ConnectedEquivalence}
Let $X$ be a limit space. $X$ is connected if and only if for every $x,y\in X$ and every local covering $\mathcal{U}$ there exists a finite number of $U_\alpha$'s in $\mathcal{U}$, namely $\{U_{\alpha_1},U_{\alpha_2},\ldots,U_{\alpha_n}\}$ such that $x\in U_{\alpha_1},y\in U_{\alpha_n}$, and $U_{\alpha_i}\cap U_{\alpha_{i+1}}\neq\emptyset$ for $1\leq i<n$.
\end{proposicion}

\begin{proof}
Let $X$ be a convergence space.
\begin{itemize}[wide]
	\item[$(\Rightarrow)$] Assume $X$ is connected. Let's take $x,y\in X$ and an arbitrary local covering $\mathcal{U}$. Let's the define the sets in $X$:
	\begin{align*}
	A_1:=&\bigcup_{U\in\mathcal{U}\mbox{ st }x\in U} U,\\ A_2:=&\bigcup_{U\in\mathcal{U}\mbox{ st }x\in U\mbox{ and }x\in A_1} U,\ldots,\\ A_n:=&\bigcup_{U\in\mathcal{U}\mbox{ st }x\in U\mbox{ and }x\in A_{n-1}} U,\ldots
	\end{align*}	
	Note that $A_1\subset A_2\subset\ldots\subset A_n\subset \ldots$ by construction. Since for every filter $\mathcal{F}\rightarrow x'$ we have that $\mathcal{F}\cap [x']\rightarrow x'$ and $\mathcal{U}$ is a local covering of $X$, then for every filter $\mathcal{F}\rightarrow x'$ we get that $A_k\in\mathcal{F}$ if $x'\in A_k$. (We can take only elements whose contains $x$ with this observation.)
	If $y\in A_n$ for some $n\in\mathbb{N}$, we can take $\{U_1,\ldots,Un\}$ with the desired conditions.
	If $y\notin A_n$ for every $n\in\mathbb{N}$, then $y\notin \cup_{n\in\mathbb{N}}A_n$. Let's take $z\in \cup_{n\in\mathbb{N}}A_n$. Thereby every $U\in\mathcal{U}$ such that $z\in U$ is inside $A_{k+1}$. Therefore $A_{k+1}\in\mathcal{F}$ for every $\mathcal{F}\rightarrow x$.
	On the other hand, let's take $z\notin\cup_{n\in\mathbb{N}}$ and $\mathcal{F}\rightarrow z$. If there exists $B\in\mathcal{F}\cap\mathcal{U}$ such that $B\cap\cup_{n\in\mathbb{N}} A_n\neq\emptyset$, there would exist $k\in\mathbb{N}$ such that $B\cap A_k\neq\emptyset$, then $z\in A_{k+1}$. Therefore every $\mathcal{F}\rightarrow z$ has an element within $(\cup_{n\in\mathbb{N}}A_n)^c$, so that set is in $\mathcal{F}$. $\Absurd$ (This means that $X$ is not connected.) Thus $y\in A_n$ for some $n\in\mathbb{N}$.
	\item[$(\Leftarrow)$] (We work with the contrapositive affirmation.) Assume $X$ is non-connected. By definition, there exist non-empty $A,B\subset X$ such that $A\cup B=X$, $A\cap B=\emptyset$ and there exists a local covering of $X$ such that every $U\in\mathcal{U}$ is a subset either $A$ or $B$
	Take $x\in A$ and $y\in B$. Let's suppose that there exists $\{U_1,\ldots,U_n\}\subset \mathcal{U}$ such that $x\in U_1$, $y\in U_n$ and $U_i\cap U_{i+1}\neq\emptyset$. Since $x\in U_1$, then $U_1\subset A$. If $U_i\subset A$, then $U_i\cap U_{i+1}\subset A$ and $U_{i+1}\subset A$. Thus, by induction, $U_n\subset A$. $\Absurd$ (This implies that $y\in A$.) Thus there does not exist any finite collection of elements of $\mathcal{U}$ with those properties for $x,y\in X$ and that local covering of $X$.\qedhere
\end{itemize}
\end{proof}

We take an analogy from \cite{Rieser_2021} to build covering limit spaces.

\begin{definicion}[Locally Trivial]
Let $p:E\rightarrow B$ be a surjective continuous function between limit spaces, and let $U$ be a subset of $B$. If $F$ is another limit space, let $q_1:U\times F\rightarrow U$ be the projection under the first factor. We say that a homeomorphism $\Phi:p^{-1}(U)\rightarrow U\times F$ is a \emph{trivilization} of $p$ over $U$ if $q_1\circ\Phi=p$. When $\Phi$ exists, $p$ is called \emph{trivial} over $U$. We say that $p$ is \emph{locally trivial} if there exists a local covering $\mathcal{U}$ of $B$ such that $p$ is trivial over $U$ for each $U\in\mathcal{U}$.
\end{definicion}

\begin{definicion}[Covering map]
\label{def:CoveringMap}
A locally trivial map $p:E\rightarrow B$ is said to be \emph{covering map} if for every $b\in B$, $F_b=p^{-1}(b)$ is discrete (that is, every point is closed and open) in the subspace limit structure on $F_b\subset E$.
\end{definicion}

\begin{proposicion}
Let $p:E\rightarrow B$ be a covering map through the local covering $\mathcal{U}$. Then $\mathcal{U}$ induces a covering map in $E$.
\end{proposicion}

\begin{proof}
Let's take $y\in E$ and $\mathcal{F}\rightarrow y$. Since $p$ is continuous, there exists $U\in p(\mathcal{F})\cap \mathcal{U}$; moreover $F_{p(y)}$ is discrete and there exists a homeomorphism $\Phi_U:p^{-1}(U)\rightarrow U\times F_{p(y)}$ such that $q_1\Phi_U=p$.

We have that $\Phi_U^{-1}(U\times\{y\})\cong_{\Phi_U} U\times \{y\}\cong_{q_1} U$ because $F_{p(y)}$ is discrete. Then $\mathcal{F}\rightarrow y$ if and only if $\Phi_U(\mathcal{F}\cap p^{-1}(U))\rightarrow (p(y),y)$ because $p^{-1}(U)\in\mathcal{F}$. This is equivalent to $\Phi_U(\mathcal{F}\cap p^{-1}(U))\cap (B\times \{y\})=p(\mathcal{F}\cap p^{-1}(U))\times \{y\}\rightarrow (p(y),y)$. Thus $\Phi_U^{-1}(U\times \{y\})\in\mathcal{F}$.

We obtain a local covering of $E$, $\{\Phi_U^{-1}(U\times \{y\}): U\in\mathcal{U}, y\in E \}$.
\end{proof}

We require a generalization of neighborhood basis to be capable to define (any) local connectedness. Thus we make the following definition.

\begin{definicion}[Local Covering System Base]
\label{def:LocalCoveringSystemBase}
Let $X$ be a limit space and $x$ a point at $X$. We say that $\mathcal{B}$ is a \emph{local covering system base} of $X$ at $x$ if for every $\mathcal{H}\rightarrow x$ there exists a $\mathcal{F}\rightarrow x$ such that $\mathcal{F}\subset \mathcal{H}$, and for every $A\in\mathcal{F}$ there exists $B\in\mathcal{B}\cap\mathcal{F}$ such that $B\subset A$.

Let $A\subset X$. By a \emph{covering system base} of $A$ is meant a collection $\mathcal{C}$ of subsets of $X$ which is a local covering system base of $X$ at each point of $A$. Whenever we omit the set, we will assume that $A=X$.
\end{definicion}

Observe that a local covering system base is in particular a local covering, however the former has the extra condition to look at within the sets in the filters, that is, to have smaller elements.

\begin{ejemplo}
\label{exam:LCS-base}
\hspace{1mm}

\begin{enumerate}
\item For every topological space, every open neighborhood base is a local covering system base. In particular, let $X=\mathbb{R}$ joint with its euclidean topology. $\mathcal{B}_t:=\{(t-\frac{1}{n},t+\frac{1}{n}):n\in\mathbb{N}\}$ is a local covering system base at $t$.

\item Let $r>0$ and $X=S^1$ joint with the limit space defined in \ref{exam:ChoquetNonCech}. Let $n\in\mathbb{N}$, then
\begin{align*}
\mathcal{B}\coloneqq \{ \{x_1,\ldots,x_k\}\mid k\geq n\text{ and }\exists x\in X\text{ such that }d(x_i,x)\leq r\ \forall i\in\{1,\ldots,k\}\}
\end{align*} 
is a local covering system base.

Observe that every filter with less than $n$ elements are including in one with $k\leq n$ elements.
\end{enumerate}
\end{ejemplo}

Since \ref{lem:Filter with point}, every element of a local covering system base of $X$ at $x$ can include $x$ by itself.

\begin{definicion}
$X$ is \emph{locally connected} if for every $x$ at $X$ there exists a local covering system base $\mathcal{B}_x$ at $x$ such that $B$ is connected for every $B\in\mathcal{B}_x$.
\end{definicion}

With the intention to achieve that every locally connected topological space is a locally connected limit space, we realize why we define local covering system base through coarser filter. Take for example $\mathbb{R}$ with its euclidean topology. The filter $[\{\{\frac{1}{k}:k\geq n\}\mid n\in\mathbb{N}\}]_\mathbb{R}\rightarrow 0$ for being finer than the neighborhood filter. However, there is not a connected space contained in $\{\frac{1}{n}:n\in\mathbb{N}\}$.

\begin{ejemplo}
\hspace{1mm}
\begin{enumerate}
\item Every locally connected topological space is a locally connected limit space.
\item Example \ref{exam:ChoquetNonCech} is a locally connected limit space. We already see in \ref{exam:LCS-base} that $\mathcal{B}$ is a local covering system base. Observe that $\{x_1,\ldots,x_k\}$ is connected because $\{\{x_1,\ldots,x_k\}\}\rightarrow x_i$ as subspace of $X$ for every $i\in\{1,\ldots,k\}$.
\end{enumerate}
\end{ejemplo}

Some of the following results are established  in \cite{Spanier_1966} for topological spaces. We keep the references changing that category for limit spaces. A few proofs only need the group structure in first homotopy group; there we do not do any change at all.

First we can characterize a covering map by covering spaces in the components of the codomain.

\begin{teorema}[\cite{Spanier_1966}; 2.1.11]
If $B$ is locally connected, a continuous map $p:E\rightarrow B$ is a covering map if and only if for each component $C$ of $B$ the map $p|p^{-1}C:p^{-1}C\rightarrow C$ is a covering map.
\end{teorema}

\begin{proof}
Let $p:E\rightarrow B$ be a continuous map
\begin{itemize}[wide]
\item[$(\Rightarrow)$] Assume $p$ is a covering map. Then there exists $\mathcal{U}$ such that $p$ is trivial over $U\in\mathcal{U}$. Let $C$ be a component of $X$ and $x\in C$. Let $\mathcal{F}$ be a filter such that $\mathcal{F}=\mathcal{F}\cap [x]$ and for every $A\in\mathcal{F}$ there exists $A'\in\mathcal{G}$ connected such that $A'\subset A$, which exists because $B$ is locally connected.

Let $U\in\mathcal{U}\cap\mathcal{F}$ and $V$ be a component of $U$ containing $x$, where $x\in U$ by construction. Since $V$ is a component in $U$, then $V\in\mathcal{F}$ and $V\subset C$. Thus $[\mathcal{F}\cap V]\rightarrow x$ on $V$ and $[\mathcal{F}\cap C]\rightarrow x$ on $C$.

Since $\mathcal{F}$ is arbitrary, all of the $V$ forms a local covering system of $C$.

\item[$(\Leftarrow)$] Assume $p|p^{-1}C: p^{-1}C\rightarrow C$ is a covering map for each component $C$ of $X$.

Let $C$ be a component, $x\in C$ and $\mathcal{F}\rightarrow x$ on $B$. Then there exists $\mathcal{G}\rightarrow x$ such that $\mathcal{G}\supset\mathcal{F}$ and for every $A\in\mathcal{G}$ there exists $A'\in\mathcal{G}$ connected such that $A'\subset A$.

Thus $C\in\mathcal{G}$, and $\mathcal{G}\cap C\rightarrow x$ on $C$. Since $p|p^{-1}C$ is a covering map, there exists $\mathcal{U}_C$ such that $p|p^{-1}C$ is trivial over every $U\in\mathcal{U}_C$. In particular, there exists $U\in\mathcal{U}_C\cap\mathcal{G}$. Because $\mathcal{G}=[\mathcal{G}\cap C]_C$, then $U\in \mathcal{G}$. Therefore, $\mathcal{U}=\bigcup\limits_{C \text{ component}}\mathcal{U}_C$ is our local covering system of $B$.\qedhere
\end{itemize}
\end{proof}

\begin{lema}[Based on \cite{Spanier_1966}; 2.1.12]
\label{lem:HomeoToComponents}
Let $p:E\rightarrow B$ a covering map through $\mathcal{U}$, and $U\in\mathcal{U}$ connected. Then $p$ maps each component of $p^{-1}(U)$ homeomorphically onto $U$.
\end{lema}

\begin{proof}
Let $x\in U$, then $U\times p^{-1}(U)$. Thus $p^{-1}(U)=\cup_{z\in p^{-1}(x)} \Phi_U^{-1}(U\times \{z\})$ where $U\times \{z\}$ is connected, therefore $\Phi^{-1}(U\times \{z\})$. Since $p$ is locally trivial over $U$, $\Phi_U^{-1}(U\times \{z\})\cong U$.
\end{proof}

\begin{corolario}[\cite{Spanier_1966}; 2.1.13]
\label{coro:SurjectiveImpliesCoveringMap}
Consider a commutative triangle
\begin{align*}
\xymatrix{
\tilde{X}_1 \ar[rr]^p \ar[rd]_{p_1} & & \tilde{X}_2 \ar[ld]^{p_2}\\
& X
}
\end{align*}
where $X$ is locally connected and $p_1$ and $p_2$ are covering maps. If $p$ is surjective, it is a covering map.
\end{corolario}

\begin{proof}
Let $\mathcal{U}_i$ be the local covering with which $p_i:\tilde{X}_i\rightarrow X$ is locally trivial. Let's take $\tilde{\mathcal{F}}\rightarrow \tilde{x}_2$, then $p_2(\tilde{\mathcal{F}})\rightarrow x$ with $x\coloneqq p_2(\tilde{x}_2)$. Therefore, there exists $\mathcal{G}\rightarrow x$ on $X$ such that $\mathcal{G}\supset p_2(\tilde{\mathcal{F}})$ and for every $A\in\mathcal{G}$ there exists $A'\in\mathcal{G}$ connected such that $A'\subset A$. Thus, there exists $U_1\in \mathcal{U}_1\cap \mathcal{G}$, $U_2\in\mathcal{U}_2\cap \mathcal{G}$ and $U\subset U_1\cap U_2\cap\mathcal{G}$ connected such that $p$ is locally trivial. In particular, $U\times p_1^{-1}(x)\cong p_1^{-1}(U)$ and $U\times p^{-1}_2(x)\cong p_2^{-1}(U)$.

Observe that, for $z\in p^{-1}(\tilde{x}_2)$ we have that $\Phi_{U,1}^{-1}(U\times \{z\})\cong_{(p_2^{-1}|U\times \{\tilde{x}_2\})\circ(p_1|U\times \{z\})}\Phi_{U,2}^{-1}(U\times \{\tilde{x}_2\})$. Therefore $p^{-1}(\Phi_{U,2}^{-1}(U\times \{\tilde{x}_2\}))\cong U\times p^{-1}(\tilde{x}_2)$. Thus $p$ is a covering map.
\end{proof}

\begin{teorema}[\cite{Spanier_1966}, 2.1.14]
If $p:\tilde{X}\rightarrow X$ is a covering map onto $X$ which is a locally connected space, then for any component $\tilde{C}$ of $\tilde{X}$ the map $p|\tilde{C}:\tilde{C}\rightarrow p(\tilde{C})$ is a covering map onto some component of $X$.
\end{teorema}

\begin{proof}
Let $\tilde{C}$ a component of $\tilde{X}$. We show the following: $p(\tilde{C})$ is connected, $p(\tilde{C})$ is a component of $X$ and $p|\tilde{C}$ is a covering map.

Since $p$ is continuous, then $p|\tilde{C}$ is continuous. Suppose that $p(\tilde{C})$ is non-continuous, then there exists disjoint subsets $A,B$ of $p(\tilde{C})$ such $\{A,B\}$ is a local covering system. Let $\tilde{\mathcal{F}}\rightarrow \tilde{x}$ on $\tilde{C}$, then $p|\tilde{C}(\tilde{F})\rightarrow x$ where $x\coloneqq p(\tilde{x})$. Therefore $A\in p|\tilde{C}(\tilde{F})$ or $B\in p|\tilde{C}(\tilde{F})$, then $p^{-1}(A)\cap \tilde{C}\in \mathcal{F}$ or $p^{-1}(B)\cap\tilde{C}\in\mathcal{F}$ $\Absurd$. This implies that $\{p^{-1}(A)\cap\tilde{C},p^{-1}(B)\cap\tilde{C}\}$ is a local covering system of $\tilde{C}$ and then $\tilde{C}$ is disconnected. Thus $p(\tilde{C})$ is connected.

Let $x\in p(\tilde{C})$ and $\mathcal{G}\rightarrow x$ such that $\mathcal{G}=\mathcal{G}\cap [x]$ and $A\in\mathcal{G}$ there exists a $A'\in\mathcal{G}$ such that $A'\subset A$. Let $B\in\mathcal{G}$ connected such that $p$ is trivial over $B$. Then $B\cap p(\tilde{C})\neq \varnothing$ and there exists $\tilde{U}$ component of $p^{-1}(B)$ such that $\tilde{U}\cap \tilde{C}\neq \varnothing$. Since $\tilde{C}$ is a component, then $\tilde{U}\subset \tilde{C}$. By \ref{lem:HomeoToComponents}, $p(\tilde{C})\subset p(\tilde{U})=B$. Observe this implies that if you add another element to $p(\tilde{C})$, then its filters contains $X-p(\tilde{C})$; thus it is disconnected.

The last argument shows that if $x\in p(\tilde{C})$ and $B\in\mathcal{G}$ connected such that $p$ is trivial over $p$, then $B\subset p(\tilde{C})$ and $(p|\tilde{C})^{-1}(B)$ is the disjoint union of those components of $p^{-1}(U)$ that meet $\tilde{C}$. It follows from \ref{lem:HomeoToComponents} that $p|\tilde{C}$ is trivial over $B$. Therefore $p|\tilde{C}$ is a covering map.
\end{proof}

We define path-connected as topology and locally path-connected. The last results in \autoref{sub:CoveringSpaces} is a set of relations between connected and path-connected spaces.

\begin{definicion}
\label{def:pathconnected}
Let $X$ be a limit space. We say that:
\begin{itemize}
	\item $X$ is \emph{path-connected} if for every $x,y\in X$ there exists $\gamma:I\rightarrow X$ continuous such that $\gamma(0)=x$ and $\gamma(1)=y$.
	\item $X$ is \emph{locally path-connected} if for every $x$ at $X$ there exists a local covering system base at $x$ $\mathcal{B}_x$ such that $B$ is path-connected for every $B\in\mathcal{B}_x$.
\end{itemize}
$C$ is a \emph{path component} if $C$ path connected and every set which contains $C$ is not path connected.
\end{definicion}

\begin{ejemplo}
\ref{exam:ChoquetNonCech} is path-connected and locally path-connected for $r>0$.

First we will observe that $\gamma:I\rightarrow S^1$ must have a finite image. Let $t\in I$ and suppose that for every $n\in\mathbb{N}$ $\#f(I\cap(t-\frac{1}{n},t+\frac{1}{n}))$ has infinite number of elements. Then for every $A\in f[I\cap (t-\frac{1}{n},t+\frac{1}{n})]$ there exists $n\in\mathbb{N}$ such that $f(I\cap (t-\frac{1}{n},t+\frac{1}{n}))\subset A$ $\Absurd$. Since $f[I\cap (t-\frac{1}{n},t+\frac{1}{n})]\rightarrow t$, then a finite subset is element of that filter. Thus there exists $n\in\mathbb{N}$ such that $f(t-\frac{1}{n},t+\frac{1}{n})$ is finite. Since $I$ is compact in topology, it is covered by a finite amount of intervals, so its image is finite.

Let $x,y\in S^1$, suppose without loss of generality that $x=0$. Let $k=\lceil \frac{y}{r} \rceil$ and $f:I\rightarrow S^1$ such that
\begin{align*}
f(t) \coloneqq \begin{cases}
0 & 0 \leq t < \frac{1}{k}\\
r & \frac{1}{k} \leq t < \frac{2}{k}\\
\vdots\\
(k-1)r & \frac{k-1}{k}\leq t < 1\\
y & t=1,
\end{cases}
\end{align*}
which is continuous by construction. Thus $S^1$ is path-connected.

To observe that $S^1$ is locally path-connected is enough to check that the elements in $\mathcal{B}$ in \ref{exam:LCS-base} are path-connected. It follows straightforward for the previous function. 
\end{ejemplo}

\begin{proposicion}
\label{prop:PathConnectedIsConnected}
Let $X$ be a path-connected limit space. Then $X$ is a connected space.
\end{proposicion}

\begin{proof}
(Proof by contradiction.) Assume that $X$ is path connected and disconnected. Since $X$ is disconnected, there exist non empty subsets $A,B$ of $X$ such that $X=A\cup B$, $A\cap B=\varnothing$ and for every filter $\mathcal{F}$ who converge to any point either $A\in\mathcal{F}$ or $B\in\mathcal{F}$. Because $A$ and $B$ are non-empty, we take $a\in A$ and $b\in B$.

$X$ is path connected, then there exists $\gamma:[0,1]\rightarrow X$ such that $\gamma(0)=a$ and $\gamma(1)=b$. Let $\mathcal{T}_t=[(t-\frac{1}{n},t+\frac{1}{n}]_I$. Since $\gamma$ is continuous, there exists $n_t\in\mathbb{N}$ such that $\gamma(t-\frac{1}{n_t},t+\frac{1}{n_t})\subset A$ or $\gamma(t-\frac{1}{n_t},t+\frac{1}{n_t})\subset B$. By topological compactness, we can take just finite many $t$'s $\Absurd$. This is a contradiction because, if we start with a $t_1$ such that $0\in (t_1-\frac{1}{n_{t_1}},1+\frac{1}{n_{t_1}})$, then $\gamma(t_1-\frac{1}{n_{t_1}},1+\frac{1}{n_{t_1}})\subset A$; thus $I\subset A$.

Thus $X$ is connected.
\end{proof}

\begin{lema}
\label{lem:ConnectedComponents}
Let $X$ be a limit space. If $X$ is locally path-connected, then $X$ has the components as path components.
\end{lema}

\begin{proof}
Let $C$ be a component. Since $X$ is locally path connected, there exists a local covering base $\mathcal{B}$ of $X$ such that each element is path connected. By \autoref{lem:ConnectedEquivalence}, between every $x,y\in C$ there exists $\{B_1,\ldots B_n\}\subset \mathcal{B}$ such that $x\in B_1$, $y\in B_n$ and $B_i\cap B_{i+1}\neq\emptyset$ for $i\in\{1,\ldots,n-1\}$. Then we build a path between $x$ and $y$. Therefore $C$ is path connected.

Observe that if there exists $C'\supset C$ path connected, then $C'$ is connected by \autoref{prop:PathConnectedIsConnected}. Thus $C'=C$ because $C$ is a component. Therefore $C$ is a path component.

On the other hand, let $C$ be a path component, then $C$ is connected by \autoref{prop:PathConnectedIsConnected}. Let $C'$ be the connected component such that $C\subset C'$. Above we show that $C'$ is a path component, then $C=C'$. Thus $C$ is a connected component.
\end{proof}

\begin{proposicion}
\label{prop:ConnLocConnPathIsPathConn}
Let $X$ a limit space. If $X$ is connected and locally path-connected, then $X$ is path-connected.
\end{proposicion}

\begin{proof}
Since $X$ is connected, $X$ has only one component. By \autoref{lem:ConnectedComponents}, that component has to be a path component. Thus $X$ is path connected.
\end{proof}

\subsection{Fibrations}

This category allows the definition of (Hurewicz) fibrations. In this section we observe some properties of these fibrations that we use later.

\begin{definicion}
\label{def:Homotopy Lifting Property}
Let $E$ and $B$ be limit spaces. A map $p:E\rightarrow B$ is said to have the \emph{homotopy lifting property} with respect to a space $X$ if, given maps $f':X\rightarrow E$ and $F':X\times I\rightarrow B$ such that $F(x,0)=pf'(x)$ for $x\in X$, there exists a map $F':X\times I\rightarrow E$ such that $F'(x,0)=f'(x)$ for $x\in X$ and $pF'=F$. It means, the following diagram commutes
\begin{align*}
\xymatrix{
X\times 0 \ar[r]^{f'} \ar[d]_{\cap} & E \ar[d]^p \\
X\times I \ar[r]_{F} \ar@{-->}[ru]^{F'} & B.
}
\end{align*}
\end{definicion}

\begin{definicion}
\label{def:fibration}
A map $p:E\rightarrow B$ is called a \emph{fibration} (or \emph{Hurewicz fiber space}) if $p$ has the homotopy lifting property with respect to every space. $E$ is called the \emph{total space} and $B$ the \emph{base space} of the fibration. For $b\in B$, $F_b:=p^{-1}(b)$ is called the \emph{fiber} of $b$.
\end{definicion}

\begin{teorema}[\cite{Spanier_1966}, 2.2.2]\label{theo:CoveringMapsAgrees}
Let $p:E\rightarrow B$ be a covering map and let $f,g:Y\rightarrow E$ liftings of the same map (that is, $pf=pg$). If $Y$ is connected and $f$ agrees with $g$ for some point of $Y$, then $f=g$.
\end{teorema}

\begin{proof}
Let $Y_1:=\{y\in Y: f(y)=g(y)\}$. Since $Y_1\neq\varnothing$, we can take $y\in Y_1$ and $\mathcal{F}\rightarrow y$ in $Y$. By continuity, $pf(\mathcal{F})\rightarrow pf(y)$ and $pg(\mathcal{F})\rightarrow pg(y)=pf(y)$. Then there exists $U\in pf(\mathcal{F})\cap pg(\mathcal{F})$ such that $p$ is trivial over $U$. Thus there exists $\tilde{U}$ such that $\tilde{U}\in g(\mathcal{F})\cap f(\mathcal{F})$, $\tilde{U}\cong_p U$ and $f(y)=g(y)\in\tilde{U}$. By construction, there exist $V,V'\in\mathcal{F}$ such that $gV=fV'=\tilde{U}$. That means that $g^{-1}(\tilde{U})\in\mathcal{F}$ and $f^{-1}(\tilde{U})\in\mathcal{F}$, what implies $g^{-1}(\tilde{U})\cap f^{-1}(\tilde{U})\in\mathcal{F}$.

If we take $z\in g^{-1}(\tilde{U})\cap f^{-1}(\tilde{U})$, then $f(z),g(z)\in \tilde{U}$. However $pf(z)=pg(z)\in U$, thus $f(z)=g(z)$. Therefore $Y_1\in\mathcal{F}$.

Let $Y_2:=\{y\in Y: f(y)\neq g(y)\}$ and $y\in Y_2$. By continuity, $pf(\mathcal{F})\rightarrow pf(y)$, $pg(\mathcal{F})\rightarrow pf(y)$ and $pf(\mathcal{F})\cap pg(\mathcal{F})\rightarrow pf(y)$. Then there exists $U\in pg(\mathcal{F})\cap pf(\mathcal{F})$ such that $p$ is trivial over $U$. Thus there exist $\tilde{U}$ and $\tilde{U}'$ such that $\tilde{U}\cong_p U \cong_p \tilde{U}'$. Therefore $g^{-1}(\tilde{U}'),f^{-1}(\tilde{U})\in\mathcal{F}$ and $g^{-1}(\tilde{U}')\cap f^{-1}(\tilde{U})\in\mathcal{F}$. We obtain that $Y_2\in\mathcal{F}$ $\Absurd$ (This contradicts that $Y$ is connected. Thus $Y_2=\varnothing$ and $Y_1=Y$.
\end{proof}

\begin{teorema}[\cite{Spanier_1966}, 2.2.3]
A covering map is a fibration.
\end{teorema}

\begin{proof}
Let $p:\tilde{X}\rightarrow X$ a covering map, and let $f':Y\rightarrow \tilde{X}$ and $F:Y\times I\rightarrow X$ be maps such that $F(y,0)=pf'(y)$, $y\in Y$.

Assume there exists a local covering $\mathcal{U}$ in $Y$ and maps $F'_U:U\times I\rightarrow \tilde{X}$ such that $F'_U(y,0)=f'(y)$, $y\in U$, and $pF'_U=F|U\times I$ for every $U\in\mathcal{U}$. Below we construct these maps and the local covering; first we show that it is defined $F'$ such that $F=pF'$ by them.

If $y\in U\cap U'$, with $U,U'\in\mathcal{U}$, then $F'_U|y\times I$ and $F'_{U'}|y\times I$ sends the path connected space $y\times I$ into $\tilde{X}$ such that for every $t\in I$
\begin{align*}
(pF'_U|y\times I)(y,t)=F(y,t)=(pF'_{U'}|y\times I)(y,t).
\end{align*}

Since $(F'_U|y\times I)(y,0)=f'(y)=(F'_{U'}|y\times I)(y,0)$ and from \ref{theo:CoveringMapsAgrees}, it follows that $F'_U|y\times I=F'_{U'}|y\times I$. Thus $F'_U|U\cap U'=F'_{U'}|U\cap U'$.

Now we define $F':Y\times I\rightarrow X$ as $F'(y,t)=F'_U(y,t)$ for $U$ such that $y\in U$. Observe $F'$ is well defined for the previous argument. We want to see $F'$ is continuous.

Let's take $\mathcal{G}\rightarrow (x,t)$, then $p_1\mathcal{G}\rightarrow x$ and $p_2\mathcal{G}\rightarrow t$ by the definition of product limit structure. Since $p_1\mathcal{G}\rightarrow x$ there exists $U\in\mathcal{U}\cap p_1\mathcal{G}$. $F'_U(\mathcal{G})\rightarrow F'_U(x,t)$ because $F'_U$ is continuous. Thus $[F'_U(\mathcal{G})]_{Y\times I}\rightarrow F'_U(x,t)=F'(x,t)$ for the definition of subspace limit structure.

If we show that $[F'_U(\mathcal{G})]_{Y\times I}= F'(\mathcal{G})$, we obtain the continuity. Observe that $U\times I\in\mathcal{G}$ and $F'_U(U\times I)=F'(U\times I)$, therefore we take elements contained in $F'(U\times I)$ without loss of generality. Let $A\subset F'(U\times I)$, then there exists $B\subset Y\times I$ such that $F'(B)=A$. Since $B\subset Y\times I$, $F'_U(B)=F'(B)$. Thus $[F'_U(\mathcal{G})]_{Y\times I}= F'(\mathcal{G})$.

Following we define those maps. Since $p$ is a covering map, there exist a local covering $\mathcal{U}$ in $X$ with which $p$ is trivial for each of its elements.

Let's take $\mathcal{F}\rightarrow y$ in $Y$ and define $\mathcal{G}_t:=[\{A\times(t-\varepsilon,t+\varepsilon) : A\in\mathcal{F},\varepsilon>0\}]$; $\mathcal{G}_t\rightarrow (y,t)$ by definition of product limit space. Therefore there exists $U_t\in\mathcal{U}\cap F(\mathcal{G}_t)$ and there exist $A_t\in\mathcal{F}$ and $\varepsilon_t>0$ such that $F(A_t\times (t-\varepsilon_t,t+\varepsilon))\subset U_t$.

By (topological) compactness of $I$ there exist $0=t_0<t_1<\ldots t_m=1$ such that $F(A'_{i}\times [t_i,t_{i+1}])\subset U'_i\in\mathcal{U}$. That condition is preserved even when we change $A'_i$ for $A:=\cap A'_i$ what also lives in $\mathcal{F}$.

Observe that $A$ depends on the filter $\mathcal{F}$, hence $A$'s could be candidates to be the local covering in $Y$. Maps are defined by induction:

\textbf{Base step:} ($G_1:A\times [0,t_1]\rightarrow \tilde{X}$) Since $p$ is trivial over $U'_1$, then $p^{-1}(U'_1)=\cup \overline{U}_j$ such that $\overline{U}_j\cong U'_1$ via $p$ restricted. Define $V_j:=f'^{-1}(\overline{U}_j)$. Observe that $pf'(\cup V_j)=U'_1=F(\cup V_j\times \{0\})$.

We claim that $A\subset \cup V_j$. Let $y\in A$, $F(y,0)=pf'(y)\in U'_1$, then there exists $j$ such that $f'(y)\in \overline{U}_j$. Thus $y\in f'^{-1}(\overline{U}_j)=V_j$.

Since $\overline{U}_j\cong_p U'_1$, we can define $G_1(y,t)=(p|V_j)^{-1}F(y,t)$ with $y\in A,t\in [0,t_1]$. By construction, $V_j\cap V_{j'}=\varnothing$, then $G_1$ is well-defined. $G_1$ is continuous because it is composition of continuous maps.

\textbf{Inductive step:} ($G_i$ from $G_{i-1}$) Assume $G_{i-1}$ is define for $2<i<m$. We know $F(A\times [t_i,t_{i+1}])\subset U'_i$. Since $p$ is trivial over $U'_i$, then $p^{-1}(U'_i)=\cup \overline{U}_j^i$ such that $\overline{U}_j^i\cong U'_i$ for every $j$. Let's define $V_k^i:=\{y\in A:G_{i-1}(y,t_i)\in \hat{U}_k^i\}$. Observe that $\cup_k V_k^i=A$. Define $G_i(y,t)=(p|V_k^i)^{-1}F(y,t)$.

Thus we build $G_i:A\times [t_i,t_{i+1}]\rightarrow \tilde{X}$, with $i\{0,\ldots,m-1\}$, such that
\begin{itemize}
	\item $pG_i=F|A\times [t_i,t_{i+1}]$.
	\item $G_i(y,0)=f'(y)$ for $y\in A$.
	\item $G_{i-1}(y,t_i)=G_i(y,t_i)$ for $y\in A$.
\end{itemize}

With this in mind, we define $F'_A:A\times I\rightarrow \tilde{X}$ such that $F'_A|A\times [t_i,t_{i+1}]=G_i$.
\end{proof}

\begin{definicion}
\label{def:unique path lifting}
A map $p:E\rightarrow B$ is said to have \emph{unique path lifting} if, given paths $\omega,\omega'$ in $E$ such that $p\omega=p\omega'$ and $\omega(0)=\omega'(0)$, then $\omega=\omega'$.
\end{definicion}

\begin{lema}[\cite{Spanier_1966}, 2.2.4]
If a map has a unique path lifting, it has the unique lifting property for path-connected spaces.
\end{lema}

\begin{proof}
Assume that $p:E\rightarrow B$ has unique path lifting. Let $Y$ be path connected and suppose that $f,g:Y\rightarrow E$ are maps such that $pf=pg$ and $f(y_0)=g(y_0)$ for some $y_0\in Y$. We must show $f=g$. Let $y\in Y$ and let $\omega$ be a path in $Y$ from $y_0$ to $y$. Then $f\omega$ and $g\omega$ are paths in $E$ that are liftings of the same path in $B$ and have the same origin. Because $p$ has unique path lifting, $f\omega=g\omega$. Therefore $f(y)=f\omega(1)=g\omega(1)=g(y)$.
\end{proof}

\begin{teorema}[\cite{Spanier_1966}, 2.2.5]
\label{theo:UniquePathLiftingNoNonConstantsPaths}
A fibration has unique path lifting if and only if every fiver has no nonconstants paths.
\end{teorema}

\begin{proof}
Assume that $p:E\rightarrow B$ is a fibration with unique path lifting. Let $\omega$ be a path in the fiber $p^{-1}(b)$ and let $\omega'$ be the constant path in $p^{-1}(b)$ such that $\omega'(0)=\omega(0)$. Then $p\omega=p\omega'$, which implies $\omega=\omega'$. Hence $\omega$ is a constant path.

Conversely, assume that $p:E\rightarrow B$ is a fibration such that every fiber has no nontrivial path and let $\omega$ and $\omega'$ be paths in $E$ such that $p\omega=p\omega'$ and $\omega(0)=\omega'(0)$. For $t\in I$, let $\omega''_t$ be the path defined by
\begin{align*}
\omega''_t(t'):=\left\lbrace \begin{array}{ll}
\omega((1-2t')t) &\ 0\leq t'\leq \frac{1}{2}\\
\omega'((2t'-1)t) &\ \frac{1}{2}\leq t' \leq 1.
\end{array}\right.
\end{align*}
Then $\omega''_t$ is a path in $E$ from $\omega(t)$ to $\omega'(t)$ and $p\omega''_t$ is a closed path in $B$ that is homotopic relative to $\{0,1\}$ to the constant path at $p\omega(t)$. By the homotopy lifting property of $p$, there is a map $F':I\times I\rightarrow E$ such that $F'(t',0)=\omega''_t(t')$ and $F'$ maps $(0\times I)\cup (I\times 1)\cup (1\times I)$ to the fiber $p^{-1}(p\omega(t))$. Because $p^{-1}(p\omega(t))$ has no nonconstant paths, $F'$ maps $0\times I$, $I\times 1$, and $1\times I$ to a single point. It follows that $F'(0,0)=F'(1,0)$. Therefore $\omega''_t(0)=\omega''_t(1)$ and $\omega(t)=\omega'(t)$.
\end{proof}

We observe that fibrations (with unique path lifting) have more structure that covering maps.

\begin{teorema}[\cite{Beattie_Butzmann_2002}, 2.2.6]
The composition of fibrations (with unique path lifting) is a fibration (with unique path lifting).
\end{teorema}

\begin{lema}[\cite{Spanier_1966}, 2.3.1]
\label{lem:FibrationRestrictedToPathComponent}
Let $p:E\rightarrow B$ be a fibration. If $A$ is any path component of $E$, then $pA$ is a path component of $B$ and $p|A:A\rightarrow pA$ is a fibration.
\end{lema}

\begin{proof}
Since $pA$ is the continuous image of a path connected space, it is path connected. It is a path component subset of $B$, as a consequence of that if $\omega$ is a path in $B$ that begins in $pA$, there is a lifting $\tilde{\omega}$ of $\omega$ that begins in $A$. Since $A$ is a path component of $E$, $\tilde{\omega}$ is a path in $A$. Therefore $\omega=p\tilde{\omega}$ is a path in $pA$.

To show that $p|A:A\rightarrow pA$ has the homotopy lifting property, let $f':Y\rightarrow A$ and $F:Y\times I\rightarrow pA$ be maps such that $F(y,0)=pf'(y)$. Because $p$ is a fibration, there is a map $F':Y\times I\rightarrow E$ such that $pF'=F$ and $F'(y,0)=f'(y)$. For any $y\in Y$, $F'$ must map $y\times I$ into the path component of $E$ containing $F'(y,0)$. Therefore $F'(y\times I)\subset A$ for all $y$, and $F':Y\times I\rightarrow A$ is a lifting of $F$ such that $F'(y,0)=f'(y)$.
\end{proof}

\begin{teorema}[\cite{Spanier_1966}, 2.3.2]
Let $p:E\rightarrow B$ be a map. If $E$ is locally path connected, $p$ is a fibration if and only if for each path component $A$ of $E$, $p(A)$ is a path component of $B$ and $p|A:A\rightarrow pA$ is a fibration.
\end{teorema}

\begin{proof}
Let $p:E\rightarrow B$ be a continuous map and $E$ be a locally path connected space.

$(\Rightarrow)$ Apply \autoref{lem:FibrationRestrictedToPathComponent} to each path component of $B$.

$(\Leftarrow)$ Let $f':Y\rightarrow E$ and $F:Y\times I\rightarrow B$ maps such that $F(y,0)=pf'(y)$ for every $y\in Y$. Let $\{A_j\}$ be the path components of $E$. Since \ref{lem:ConnectedComponents}, $\{A_j\}$ are disjoint components because $E$ is locally path connected. Denote $V_j:=f'^{-1}(A_j)$ and observe that $\cup V_j=Y$. Then it suffices to construct maps $F'_j:V_j\times I\rightarrow E$ for all $j$ such that $pF'_j=F|V_j\times I$ and $F_j(y,0)=f'(y)$.

Since $y\times I$ is path connected, $F(y\times I)$ is contained in the path component of $B$ containing $F(y,0)=pf'(y)$. Because $pA_j$ is a path component of $B$, then $F(V_j\times I)=F(\cup_{y\in V_j}y\times I)\subset pA_j$ for each $j$.

Because $p|A_j:A_j\rightarrow pA_j$ is a fibration, there exists $F_j':V_j\times I\rightarrow A_j$ such that $pF'_j:F|V_j\times I$ and $F'_j(y,0)=f'(y)$ for $y\in V_j$. Therefore $p$ has the homotopy lifting property.
\end{proof}

We can add ``with unique path lifting'' to the theorem. This is possible because every path lies in some path component of the space.

\begin{lema}[\cite{Spanier_1966}, 2.3.3]
\label{lem:FibrationWithUniquePathLiftingPreImageHomotopy}
Let $p:\tilde{X}\rightarrow X$ be a fibration with unique path lifting. If $\omega$ and $\omega'$ are paths in $\tilde{X}$ such that $\omega(0)=\omega'(0)$ and $p\omega\simeq p\omega'$, then $\omega\simeq\omega'$.
\end{lema}

\begin{proof}
Let $F:I\times I\rightarrow X$ be a homotopy relative to $\{0,1\}$ from $p\omega$ to $p\omega'$. By the homotopy lifting property of fibrations, there is a map $F':I\times I\rightarrow \tilde{X}$ such that $F'(t,0)=\omega(t)$ and $pF'=F$. Then $F'(0\times I)$ and $F(1\times I)$ are contained in $p^{-1}(p\omega(0))$ and $p^{-1}(p\omega(1))$ respectively. By \autoref{theo:UniquePathLiftingNoNonConstantsPaths}, $F'(0\times I)$ and $F'(1\times I)$ are single points. Hence $F'$ is a homotopy relative to $\{0,1\}$ from $\omega$ to some path $\omega''$ such that $\omega''(0)=\omega(0)$ and $p\omega''=p\omega'$. Since $\omega'(0)=\omega(0)$, it follows from the unique-path-lifting property of $p$ that $\omega'=\omega''$ and $F':\omega\simeq \omega'\ rel\ \{0,1\}$.
\end{proof}

\begin{teorema}[\cite{Spanier_1966}, 2.3.4]
Let $p:\tilde{X}\rightarrow X$ be a fibration with unique path lifting. For any $\tilde{x}_0\in E$ the homomorphism
\begin{align*}
p_\#:\pi(\tilde{X},\tilde{x}_0)\rightarrow \pi(X,x)
\end{align*}
is a monomorphism.
\end{teorema}

\begin{proof}
Apply \autoref{lem:FibrationWithUniquePathLiftingPreImageHomotopy} to closed paths $\omega,\omega'$ in $\tilde{X}$ with origin $\tilde{x}_0$ such that $[p\omega]=[p\omega']$. Thus $[\omega]=[\omega']$.
\end{proof}

In \cite{Spanier_1966} page 45, we have the following definition of $h_{[\omega]}([\tau])=[\omega * \tau * \omega^{-1}]$ where $\tau$ is closed path in the end of $\omega$.

\begin{lema}[\cite{Spanier_1966}, 2.3.5]
Let $p:E\rightarrow B$ be a fibration with unique path lifting and assume that $E$ is a nonempty path-connected space. If $e_0,e_1\in E$, there is a path $\omega$ in $B$ from $p(e_1)$ to $p(e_1)$ such that
\begin{align*}
p_\#\pi(E,e_0)=h_{[\omega]}p_\#\pi(E,e_1).
\end{align*}
Conversely, given a path $\omega$ in $B$ from $p(e_0)$ to $x_1$, there is a point $e_1\in p^{-1}(x_1)$ such that
\begin{align*}
h_{[\omega]}p_\#\pi(E,e_1)=p_\#\pi(E,e_0).
\end{align*}
\end{lema}

\begin{proof}
For the first part, let $\tilde{\omega}$ be a path in $E$ from $\tilde{x}_0$ to $\tilde{x}_1$. Then $\pi(E,\tilde{x}_0)=h_{[\tilde{\omega}]}\pi(E,\tilde{x}_1)$. Therefore $p_\#\pi(E,\tilde{x}_0)=h_{[p\tilde{\omega}]}p_\#\pi(E,\tilde{x}_1)$ and so $p\tilde{\omega}$ will do as the path from $p(\tilde{x}_0)$ to $p(\tilde{x}_1)$.

Conversely, given a path $\omega$ in $B$ from $p(\tilde{x}_0)$ to $x_1$, let $\tilde{\omega}$ be a path in $E$ such that $\tilde{\omega}(0)=\tilde{x}_0$ and $p\tilde{\omega}=\omega$. If $\tilde{x}_1=\tilde{\omega}(1)$, then $h_{[\omega]}p_\#\pi(E,\tilde{x}_1)=p_\#(h_{[\tilde{\omega}]}\pi(E,\tilde{x}_1))=p_\#\pi(E,\tilde{x}_0)$.
\end{proof}

The following result is direct from the last lemma.

\begin{teorema}[\cite{Spanier_1966}, 2.3.6]
\label{theo:FibrationwUPLconjugacyClass}
Let $p:E\rightarrow B$ be a fibration with unique path lifting and assume that $E$ is a nonempty path-connected space. For $x_0\in pE$ the collection $\{p_\#\pi(E,e_0):e_0\in p^{-1}(x_0)\}$ is a conjugacy class in $\pi(B,x_0)$. If $\omega$ is a path in $pE$ from $x_0$ to $x_1$, then $h_{[\omega]}$ maps the conjugacy class in $\pi(B,x_1)$ to the conjugacy class in $\pi(B,x_0)$.
\end{teorema}

\subsection{Lifting Theorem}

\begin{lema}[\cite{Spanier_1966}, 2.4.1]
Let $p:E\rightarrow B$ a fibration. Any map of a contractible space to $B$ whose image is contained in $p(E)$ can be lifted to $E$.
\end{lema}

\begin{proof}
Let $Y$ be contractible and let $f:Y\rightarrow B$ be a map such that $f(Y)\subset p(E)$. Because $Y$ is contractible, $f$ is homotopic to a constant map of $Y$ to some point of $f(Y)$. $f(Y)\subset p(E)$, so this constant map cam be lifted to $E$. The homotopy lifting property then implies that $f$ can be lifted to $E$.
\end{proof}

\begin{lema}[\cite{Spanier_1966}, 2.4.2]
\label{lem:FibrationWithUniquePathLiftingStrongDeformation}
Let $p:(\tilde{X},\tilde{x}_0)\rightarrow (X,x_0)$ be a fibration with unique path lifting. If $y_0$ is a strong deformation retract of $Y$, any map $(Y,y_0)\rightarrow (X,x_0)$ can be lifted to a map $(Y,y_0)\rightarrow (\tilde{X},\tilde{x}_0)$.
\end{lema}

\begin{proof}
Let $f:(Y,y_0)\rightarrow (X,x_0)$ be a map. $f$ is homotopy relative to $y_0$ to the constant map $y\mapsto x_0$. The constant map can be lifted to the constant map $y\mapsto\tilde{x}_0$. By the homotopy lifting property, $f$ can be lifted to a map $f':Y\rightarrow \tilde{X}$ such that $f'$ is homotopic to the constant map $y\mapsto \tilde{x}_0$ by a homotopy which maps $y_0\times I$ to $p^{-1}(x_0)$. Because $p^{-1}(x_0)$ has no nonconstant maps by \autoref{theo:UniquePathLiftingNoNonConstantsPaths}, $f'(y_0)=\tilde{x}_0$.
\end{proof}

\begin{lema}[\cite{Spanier_1966}, 2.4.3]
\label{lem:StrongDeformationRetractPYy0}
The constant path at $y_0$ is a strong deformation retract of the path space $P(Y,y_0)$.
\end{lema}

\begin{proof}
A strong deformation retraction $F:P(Y,y_0)\times I\rightarrow P(Y,y_0)$ to the constant path at $y_0$ is defined by $F(\omega,t)(t')=\omega((1-t)t')$ where $\omega\in P(Y,y_0)$ and $t,t'\in I$.
\end{proof}

Consider $P(Y,y_0)$ as \autoref{def:PathSpace}. We define $\varphi:P(Y,y_0)\rightarrow Y$ as $\omega\mapsto \omega(1)$. $\varphi$ is surjective if and only if $Y$ is path connected. To observe continuity, take $\mathcal{H}\rightarrow f$ in $P(Y,y_0)$. By continuous limit structure definition, $\omega_{I,Y}(\mathcal{H}\times [1])\rightarrow f(1)$. Then $\varphi\mathcal{H}\rightarrow f(1)$. Thus $\varphi$ is continuous.

\begin{teorema}[\cite{Spanier_1966}, 2.4.4]
\label{theo:PathSpaceIsQuotient}
A connected locally path-connected space $Y$ is the quotient space of its path space $P(Y,y_0)$ by the map $\varphi$.
\end{teorema}

\begin{proof}
Let $\lambda$ be the limit structure of $Y$. We need to show that $\lambda=(\lambda_c)_q$, where $(\lambda_c)_q$ is the quotient limit structure with respect to $\varphi$.

Since $(\lambda_c)_q$ is the final structure on $Y$ of $\varphi$ and $\lambda$ is continuous, then $(\lambda_c)_q$ is finer than $\lambda$.

We prove the other contention. Let $\mathcal{F}\rightarrow x$ under $\lambda$. Without loss of generality, allow that $\mathcal{F}\cap [x]=\mathcal{F}$ and for every $A\in\mathcal{F}$ there exists $A'\in\mathcal{F}$ path-connected such that $A'\subset A$. Since $X$ is path connected, there exists $\gamma\in P(Y,y_0)$ such that $\gamma(1)=x$. Because $Y$ is locally path connected, there exists a local covering base $\mathcal{B}_x\subset\mathcal{F}$ at $x$ such that $B\in\mathcal{B}_x$ is path-connected.

Let $B\in\mathcal{B}_x$, we define $B^\star:=\{\tau:\exists \tau', \tau'(0)=x,\tau'(1)\in B$, $\tau'(0,1]\subset B,\tau=\gamma\star\tau'\}$ and $\mathcal{B}^\star:=\{B^\star: B\in\mathcal{B}_x\}$. Observe that $\gamma\in B^\star$ for every $B$.

Let $B_1,B_2\in\mathcal{B}_x$, then there exists $B_3\in\mathcal{B}$ such that $B_3\subset B_1\cap B_3$. Thus $B_1^\star \cap B_2^\star$ is nonempty because contains $B_3^\star$.

By construction, $\varphi\mathcal{B}^\star=\mathcal{B}_x$ and $\mathcal{F}=[\mathcal{B}_x]$. It just remains to prove that $[\mathcal{B}^\star]\in\lambda_c$. (Observe that $[\mathcal{B}^\star]$ would converge to $\gamma * \varepsilon_x$, where $\varepsilon_x$ is the constant path of $x$.)

Let's take the coarser filter what converges to $t$ in I$, \mathcal{T}_t:=[(t-\frac{1}{n},t+\frac{1}{n})]_I$. If $t\in [0,\frac{1}{2})$, we get that $\tau(t)=\gamma(2t)$ for every $\tau\in B^\star\in \mathcal{B}^\star$. Thus, by continuity of $\gamma$, $\omega_{I,X}([\mathcal{B}^\star\times\mathcal{T}_t])\rightarrow \gamma(2t)$.
We now check that $\omega([\mathcal{B}^\star\times\mathcal{T}_t])\rightarrow x$ for $t\in [\frac{1}{2},t]$. We divide it in two case.
\begin{itemize}
	\item $t=\frac{1}{2}$. We divide our filter in two: $[(t-\frac{1}{n},t]]$ and $[(t,t+\frac{1}{n})]$ By construction $\omega_{I,X}([\mathcal{B}^\star\times (t-\frac{1}{n}]])=[\gamma(2t-\frac{2}{n},2t]]\rightarrow \gamma(1)$.
	On the other hand, take $B^\star\in\mathcal{B}^\star$. By definition, $\omega_{I,X}(B^\star\times(t,t+1/n))\subset B$. Then $\omega_{I,X}([\mathcal{B}^\star\times (t,t+\frac{1}{n})])\supset [\mathcal{B}]=\mathcal{F}\rightarrow x$.
	\item $t\in (\frac{1}{2},1]$. For some large enough $n$, we obtain that $\omega_{I,X}(B^\star\times (t-\frac{1}{n},t+\frac{1}{n}))\subset B$. Thus the argument is analogue to the last one.
\end{itemize}
Hence $[\mathcal{B}^\star]\in\lambda_c(\gamma\star\epsilon_x)$ and $\varphi([\mathcal{B}^\star])=[\mathcal{B}_x]=\mathcal{F}\in (\lambda_c)_q(x)$. Therefore, we have that $\lambda=(\lambda_c)_q$.
\end{proof}

We highlight the following part of the proof \ref{theo:PathSpaceIsQuotient} because it implies that it is enough to prove that $\varphi$ is continuous to observe that $Y$ is the final structure  with respect to $\varphi$, that is, it is the quotient structure on $Y$.

\begin{proposicion}
\label{prop:AlmostQuotient}
Let $Y$ be a connected locally path-connected limit space, let $Z$ be a limit space. If $\overline{f}:P(Y,y_0)\rightarrow Z$ is continuous and there exists $f'$ such that $f'\varphi=\overline{f}$, then $f'$ is continuous.
\end{proposicion}

\begin{proof}
Let $\mathcal{F}\rightarrow x$ in $X$. By the result above, there exists a filter $\mathcal{H}\rightarrow\gamma$ such that $\gamma(1)=x$ and $\varphi(\mathcal{H})=\mathcal{F}$. Since $\overline{f}$ is continuous, then $\overline{f}(\mathcal{H})\rightarrow \overline{f}(\gamma)$, what is equivalent to $f'\varphi(\mathcal{H})=f'(\mathcal{F})\rightarrow f'(x)=f'\varphi(\gamma)$.
\end{proof}

\begin{teorema}[Lifting theorem; \cite{Spanier_1966}, 2.4.5]
\label{theo:LiftingTheorem}
Let $p:(\tilde{X},\tilde{x}_0)\rightarrow (X,x_0)$ be a fibration with unique path lifting. Let $Y$ be a connected locally path-connected space. A necessary. and sufficient condition that a map $f:(Y,y_0)\rightarrow (X,x_0)$ have a lifting $(Y,y_0)\rightarrow (\tilde{X},\tilde{x}_0)$ is that in $\pi(X,x_0)$
\begin{align*}
f_\#\pi(Y,y_0)\subset p_\#\pi(\tilde{X},\tilde{x}_0).
\end{align*}
\end{teorema}

\begin{proof}
If $f':(Y,y_0)\rightarrow (\tilde{X},\tilde{x}_0)$ is a lifting of $f$, then $f=pf'$ and $f_\#\pi(Y,y_0)=p_\#f'_\#\pi(Y,y_0)\subset p_\#\pi(\tilde{X},\tilde{x}_0)$, which shows that the condition is necessary.

We now prove that the condition is sufficient. If follows from \autoref{lem:StrongDeformationRetractPYy0} and \autoref{lem:FibrationWithUniquePathLiftingStrongDeformation} that if $\omega_0$ is the constant path at $y_0$, the composite
\begin{align*}
\xymatrix{
(P(Y,y_0),\omega_0)\ar[r]^\varphi & (Y,y_0) \ar[r]^f & (X,x_0)
}
\end{align*}
can be lifted to a map $\tilde{f}:(P(Y,y_0),\omega_0)\rightarrow (\tilde{X},\tilde{x}_0)$. We show that if $f_\#\pi(Y,y_0)\subset p_\#\pi(\tilde{X},\tilde{x}_0)$ and if $\omega,\omega'\in P(Y,y_0)$ are such that $\varphi(\omega)=\varphi(\omega')$, then $\tilde{f}(\omega)=\tilde{f}(\omega')$. Let $\overline{\omega}$ and $\overline{\omega}'$ be that path in $P(Y,y_0)$ from $\omega_0$ to $\omega$ and $\omega'$, respectively, defined by $\overline{\omega}(t)(t')=\omega(tt')$ and $\overline{\omega}'(t)(t')=\omega'(tt')$. Then $\tilde{f}\overline{\omega}$ and $\tilde{f}\overline{\omega}'$ are paths in $\tilde{X}$ from $\tilde{x}_0$ to $\tilde{f}(\omega)$ and $\tilde{f}(\omega')$, respectively, such that $p\tilde{f}\overline{\omega}=f\varphi\overline{\omega}=f\omega$ and $p\tilde{f}\overline{\omega}'=f\omega'$.

Because $\omega*\omega'^{-1}$ is a closed path in $Y$ at $y_0$ and $f_\#\pi(Y,y_0)\subset p_\#(\tilde{X},\tilde{x}_0)$, there is a closed path $\tilde{\omega}$ in $\tilde{X}$ at $\tilde{x}_0$ such that $(f\omega)*(f\omega')^{-1}\simeq p\tilde{\omega}$. Then $p(\tilde{f}\overline{\omega})=f\omega\simeq (p\tilde{\omega})*(f\omega')=p(\tilde{\omega}*(\tilde{f}\tilde{\omega}'))$.

By \autoref{lem:FibrationWithUniquePathLiftingPreImageHomotopy}, $\tilde{f}\overline{\omega}\simeq \tilde{\omega}*(\tilde{f}\tilde{\omega}')$. In particular, the endpoint of $\tilde{f}\overline{\omega}$, which is $\tilde{f}(\omega)$, equals to the endpoint of $\tilde{f}\overline{\omega}'$, which is $\tilde{f}(\omega')$.

It follows that there is a function $f':(Y,y_0)\rightarrow (\tilde{X},\tilde{x}_0)$ such that $f'\varphi=\tilde{f}$, and using \autoref{theo:PathSpaceIsQuotient}, we see that $f'$ is continuous. Because $pf'\varphi=p\tilde{f}=f\varphi$ and $\varphi$ is surjective, $pf'=f$. Therefore $f'$ is a lifting of $f$.
\end{proof}

\section{Universal Covering Space}
\label{sec: Universal covering space}

In the last section about limit spaces, we introduce universal covering spaces, its main homotopy property and its construction for connected locally path-connected limit spaces with an extra condition.

\subsection{Definition and properties}

Let $X$ be a connected space. The \emph{category of connected covering spaces of $X$} \cite{Spanier_1966} has objects which are covering maps $p:\tilde{X}\rightarrow X$, where $\tilde{X}$ is connected, and morphisms which are commutative triangles

\begin{align*}
\xymatrix{
\tilde{X}_1 \ar[rr]^f \ar[rd]_{p_1} & & \tilde{X}_2 \ar[ld]^{p_2}\\
& X
}
\end{align*}

Observe that if $X$ is in addition locally path-connected, then $\tilde{X}$ is also locally path-connected since $p$ is a covering map.

\begin{lema}[\cite{Spanier_1966}, 2.5.1]
\label{lem:MorphismIsACoveringMap}
In the category of connected covering spaces of a connected locally path-connected space every morphism is itself a covering projection.
\end{lema}

\begin{proof}
Consider a commutative triangle
\begin{align*}
\xymatrix{
\tilde{X}_1 \ar[rr]^f \ar[rd]_{p_1} & & \tilde{X}_2 \ar[ld]^{p_2}\\
& X
}
\end{align*}
where $p_1$ and $p_2$ are covering projections and $X$ is locally path connected. If we prove that $f$ is surjective, then $f$ is a covering projection from \ref{coro:SurjectiveImpliesCoveringMap}.
Since $\tilde{X}_2$ is path connected for \ref{prop:ConnLocConnPathIsPathConn}. Let $\tilde{x}_1\in\tilde{X}_1$ and $\tilde{x}_2\in\tilde{X}_2$ be arbitrary and let $\tilde{\omega}_2$ be a path in $\tilde{X}_2$ from $f(\tilde{x}_1)$ to $\tilde{x}_2$. Because $p_1$ is a fibration, there is a path $\tilde{\omega}_1$ in $\tilde{X}_1$ beginning at $\tilde{x}_1$ such that $p_1\tilde{\omega}_1=p_2\tilde{\omega}_2$. By the unique path lifting of $p_2,$ $f\tilde{\omega}_1=\tilde{\omega}_2$. Therefore
\begin{align*}
f(\tilde{\omega}_1(1))=\tilde{\omega}_2(1)=\tilde{x}_2
\end{align*} 
proving that $f$ is surjective.
\end{proof}

\begin{teorema}[\cite{Spanier_1966}, 2.5.2]
\label{theo:Spanier252}
Let $p_1:\tilde{X}_1\rightarrow X$ and $p_2:\tilde{X}_2\rightarrow X$ be objects in the category of connected covering spaces of a connected locally path connected space $X$. The following are equivalent:
\begin{enumerate}[label=\alph*)]
	\item There is a covering map $f:\tilde{X}_1\rightarrow \tilde{X}_2$ such that $p_2f=p_1$.
	\item For every $\tilde{x}_1\in\tilde{X}_1$ and $\tilde{x}_2\in\tilde{X}_2$ such that $p_1(\tilde{x}_1)=p_2(\tilde{x}_2)$, $p_{1\#}\pi(\tilde{X}_1,\tilde{x}_1)$ is conjugate in $\pi(X,p_1(\tilde{x}_1))$ to a subgroup of $p_{2\#}\pi(\tilde{X}_2,\tilde{x}_2)$.
	\item There exist $\tilde{x}_1\in\tilde{X}_1$ and $\tilde{x}_2\in\tilde{X}_2$ such that $p_1(\tilde{x}_1)=p_2(\tilde{x}_2)$, $p_{1\#}\pi(\tilde{X}_1,\tilde{x}_1)$ is conjugate in $\pi(X,p_1(\tilde{x}_1))$ to a subgroup of $p_{2\#}\pi(\tilde{X}_2,\tilde{x}_2)$.
\end{enumerate}
\end{teorema}

\begin{proof}
\textbf{a)} $\Rightarrow$ \textbf{b)} Given $f:\tilde{X}_1\rightarrow \tilde{X}_2$ such that $p_2f=p_1$, if $\tilde{x}_1\in\tilde{X}_1$ and $\tilde{x}_2\in\tilde{X}_2$ are such that $p_1(\tilde{x}_2)=p_2(\tilde{x}_2)$ then $p_{1\#}\pi(\tilde{X}_1,\tilde{x}_1)=p_{2\#}f_\#\pi(\tilde{X}_1,\tilde{x}_1)\subset p_{2\#}\pi(\tilde{X}_2,\tilde{x}_2)$.

Because $f(\tilde{x}_1)$ and $\tilde{x}_2$ lie in the same fiber of $p_2:\tilde{X}_2\rightarrow X$, it follow from \autoref{theo:FibrationwUPLconjugacyClass} that $p_{2\#}\pi(\tilde{X}_2,f(\tilde{x}_1))$ and $p_{2\#}\pi(\tilde{X}_2,\tilde{x}_2)$ are conjugate in $\pi(\tilde{X}_1,p_1(\tilde{x}_1))$.

\textbf{b)} $\Rightarrow$ \textbf{c)} The condition is fulfilled by every pair $\tilde{x}_1\in\tilde{X}_1$ and $\tilde{x}_2\in\tilde{X}_2$, in particular there exists one pair which fulfill the condition.

\textbf{c)} $\Rightarrow$ \textbf{a)} Assume that $\tilde{x}_1\in\tilde{X}_1$ and $\tilde{x}_2\in\tilde{X}_2$ are such that $p_1(\tilde{x}_1)=p_2(\tilde{x}_2)$ and $p_{1\#}\pi(\tilde{X}_1,\tilde{x}_1)$ is conjugate in $\pi(X,p_1(\tilde{x}_1)$ to a subgroup of $p_{2\#}\pi(\tilde{X}_2,\tilde{x}_2)$.

By \autoref{theo:FibrationwUPLconjugacyClass}, there is a point $\tilde{x}_2'\in\tilde{X}_2$ such that $\tilde{x}_2'\in p_2^{-1}(\tilde{x}_2)$ and such that $p_{1\#}\pi(\tilde{X}_1,\tilde{x}_1)\subset p_{2\#}\pi(\tilde{X}_2,\tilde{x}'_2)$.

Because $\tilde{X}_1$ is a connected locally path connected space, the \autoref{theo:LiftingTheorem} implies the existence of a map $f:(\tilde{X}_1,\tilde{x}_1)\rightarrow (\tilde{X}_2,\tilde{x}_2)$ such that $p_2f=p_1$.
\end{proof}

\begin{corolario}[\cite{Spanier_1966}, 2.5.3]
\label{coro:Spanier253}
Two objects in the category of connected covering spaces of a connected locally path connected space $X$ are equivalent if and only if their fundamental groups map to conjugate subgroups of the fundamental group of $X$.
\end{corolario}

\begin{proof}
Let $p_1:\tilde{X}_1\rightarrow X$ and $p_2:\tilde{X}_2\rightarrow X$ be objects in the category of connected covering spaces of a connected locally path connected space $X$

For the first side, let's assume they are equivalent. Then there exist $f:\tilde{X}_1\rightarrow\tilde{X}_2$ isomorphism such that $p_1=fp_2$. By \autoref{theo:Spanier252}, there exists $\tilde{x}_1\in\tilde{X}_1$ and $\tilde{x}_2\in \tilde{X}_2$ such that $p_1(\tilde{x}_1)=p_2(\tilde{x}_2)=p_2f(\tilde{x}_1)$, $p_{1\#}\pi(\tilde{X}_1,\tilde{x}_1)$ is conjugate in $\pi(X,p_1(\tilde{x}_1))$ to a subgroup of $p_{2\#}\pi(\tilde{X}_2,\tilde{x}_2)$.

On the other hand, we satisfies b) in \autoref{theo:Spanier252}, then there exists $f:\tilde{X}_1\rightarrow \tilde{X}_2$ such that $p_2f=p_1$. Moreover, $f$ and $\Id_{\tilde{X}_2}$ satisfies the same, so there exists $g$ such that $\Id_{\tilde{X}_2}=fg$. However, $g$ and $\Id_{\tilde{X}_1}$ also satisfies b), then there exists $h$ such that $\Id_{\tilde{X}_1}=gh$. Thus $h=f$, and the objects are equivalent.
\end{proof}

\begin{definicion}[\cite{Spanier_1966}]
A \emph{universal covering space} of a connected space $X$ is an object of the category of connected covering spaces of $X$ such that for any object $p':\tilde{X}'\rightarrow X$ of this category there is a morphism
\begin{align*}
\xymatrix{
\tilde{X} \ar[rr]^f \ar[rd]_p && \tilde{X}' \ar[ld]^{p'}\\
& X
}
\end{align*}
in the category.
\end{definicion}

\begin{corolario}[\cite{Spanier_1966}, 2.5.6]
Two universal covering spaces of a connected locally path connected space are equivalent.
\end{corolario}

\begin{proof}
Let $p_1:\tilde{X}_1\rightarrow X$ and $p_2:\tilde{X}_2\rightarrow X$ universal covering spaces, then there exists $f:\tilde{X}_1\rightarrow \tilde{X}_2$ and $g:\tilde{X}_2\rightarrow\tilde{X}_1$ such that $p_1=p_2f$ and $p_2=p_1g$ by definition.

For \ref{theo:Spanier252}, their fundamental groups map to conjugate subgroups of the fundamental group of $X$. Thus, both objects are equivalents in the category of connected covering spaces of $X$ since \ref{coro:Spanier253}
\end{proof}

\begin{definicion}
\label{def:simplyconnected}
Let $X$ be a path connected limit space. We say that $X$ is \emph{simply connected} if $\pi(X)=0$.
\end{definicion}

\begin{corolario}[\cite{Spanier_1966}, 2.5.7]
A simply connected space of a connected locally path connected space $X$ is a universal cover space of $X$.
\end{corolario}

\begin{proof}
It follows from $p_{1\#}\pi(\tilde{X}_1,\tilde{x}_1)\cong 0$.
\end{proof}

\subsection{Construction}

Finally we provide the definition of the universal covering space for a limit space. We need a last condition for this construction: there exists a local covering system where the loops of its elements are contractible within all of the space.

\begin{definicion}
Let $X$ be a limit space. $X$ is said to be \emph{semi-locally simply connected} if there exists a local covering system $\mathcal{U}$ such that for every $U\in\mathcal{U}$ each continuous function $\gamma:I\rightarrow U$ such that $\gamma(0)=\gamma(1)$ is contractible to the constant path at $\gamma(0)$ within $X$.
\end{definicion}

\begin{ejemplo}
\hspace{1mm}
\begin{enumerate}
\item Every semi-locally simply connected topological space is semi-locally simply connected limit space.
\item $X=S^1$ from \ref{exam:ChoquetNonCech} is semi-locally simply connected. Provide the $r>0$, the local covering system $\{E_x\mid x\in X\}$ where $E_x=\{y\in X\mid d(y,x)\leq \frac{r}{2}\}$ and every loop within $E_x$ is contractible within $X$.
\end{enumerate}
\end{ejemplo}

For the rest of the section, we assume that $X$ is a connected locally path connected semi-locally simply connected limit space. The construction of the universal covering space is given through the following space:

\begin{definicion}
Let $X$ be a limit space and $x_0\in X$. We define $\mathcal{C}:=\{[\gamma]:\gamma\in P(X,x_0)\}$, where $\gamma'\in [\gamma]$ if and only if $\gamma'(1)=\gamma(1)$ and $\gamma\simeq \gamma'$, with quotient limit structure given by $\gamma\mapsto [\gamma]$.
\end{definicion}

\begin{lema}
$X$ is the quotient space of $\mathcal{C}$ by the map $\bar{\varphi}$ such that $\bar{\varphi}[\gamma]=\gamma(1)$.
\end{lema}

\begin{proof}
Let's call $q:P(X,x_0)\rightarrow\mathcal{C}$ is the quotient map. Observe that $\bar{\varphi}$ is well-defined by definition, that is, all of the elements in the equivalence relation have the same end-point.

Since $\gamma\mapsto [\gamma]$ is quotient, we obtain that $\bar{\varphi}$ is continuous because $\varphi=\bar{\varphi}q$ and $\varphi$ is continuous, \ref{theo:PathSpaceIsQuotient}.

Let $\mathcal{F}\rightarrow x$ under $\lambda$, what is the limit structure of $X$. Without loss of generality, $\mathcal{F}\cap[x]=\mathcal{F}$ and for every $A\in\mathcal{F}$ there exists $A'\in\mathcal{F}$ path-connected such that $A'\subset A$. There exists $\gamma\in P(X,x_0)$ such that $\gamma(1)=x$ because $X$ is path-connected. Since $X$ is locally path connected and semi-locally simply connected, then there exists a local covering base at $x$ $\mathcal{B}\subset\mathcal{F}$ such that $[\mathbb{B}]_X=\mathcal{F}$, and $B\in\mathcal{B}$ is path connected and every closed path in $B$ is contractible within $X$.

Let's define $B^\star:=\{\tau:\exists \tau',\tau'(0)=x,\tau'\in B,\tau'(0,1]\subset B,\tau=\gamma*\tau'\}$ for every $B\in\mathcal{B}$ and $\mathcal{B}^\star:=\{B^\star:B\in\mathcal{B}\}$. By \autoref{theo:PathSpaceIsQuotient}, $[\mathcal{B}^\star]_{P(X,x_0)}\rightarrow \gamma*\varepsilon_x$. Then the quotient under homotopy, denoted by $B^\star_{[\ ]}$, converges to $[\gamma]$ and their image under $\bar{\varphi}$ converges to $\gamma(1)$. Hence, $\bar{\varphi}$ is a quotient map.
\end{proof}

Once again, we have the following proposition as a corolary of the last lemma. The proof is completely analogous to \ref{prop:AlmostQuotient}.

\begin{proposicion}
Let $X$ be a connected locally path-connected semi-locally simply connected space, let $Z$ be a limit space and $\bar{\varphi}:\mathcal{C}\rightarrow X$ such that $\tau\mapsto \tau(1)$. If $\bar{f}:\mathcal{C}\rightarrow Z$ is continuous and there exists $f'$ such that $f'\bar{\varphi}=\bar{f}$, then $f'$ is continuous.
\end{proposicion}

\begin{lema}
Let $x\in X$, then $\bar{\varphi}^{-1}(x)$ is discrete under subspace limit structure.
\end{lema}

\begin{proof}
By definition, $\bar{\varphi}^{-1}(x):=\{\tau\in \mathcal{C}: \tau(1)=x\}$. Let $\mathcal{F}\rightarrow[\tau]$ in $\bar{\varphi}^{-1}(x)$, and then $[\mathcal{F}]_\mathcal{C}\rightarrow [\tau]$ in $\mathcal{C}$. Since $X$ is locally path connected and semi-locally simply connected, then there exists $\mathcal{G}$ in $X$ such that $[\mathcal{F}]_\mathcal{C}=\mathcal{G}^\star$ (this argument is the same to the one in the proof of \autoref{Teo:HomFilterMod}) and a subset $A\subset X$ containing $x_0$ and with every loop (closed path) contained there contractible. Thus, if we observe only the closed path in $\mathcal{G}^\star$ all of them are homotopy to one element, say $\tau$, and the filter is the filter generated just by $\tau$.
\end{proof}

\begin{lema}
\label{lem:VarphiIsLocallyTrivialMap}
$\bar{\varphi}$ is locally trivial map.
\end{lema}

\begin{proof}
Since $X$ is semi-locally simply connected space, there exists a local covering $\mathcal{U}$ such that for every $U\in\mathcal{U}$ every closed path at $U$ is contractible in $X$.

Let $x\in X$ and, without loss of generality, let $\mathcal{F}\rightarrow x$ such that $\mathcal{F}=\mathcal{F}\cap[x]$ and for every $A\in\mathcal{F}$ there exists $A'\in\mathcal{F}$ path-connected such that $A'\subset A$. Then there exists $U(x,\mathcal{F})\in\mathcal{U}\cap\mathcal{F}$ and $B(x,\mathcal{F})\in \mathcal{F}$ such that $B(x,\mathcal{F})\subset U(x,\mathcal{F})$. Note that $B(x,\mathcal{F})$ inherits that every cycle in $B$ is contractible in $X$. By simplicity, we write $B(x,\mathcal{F})$ as $B$, having in mind that $B$ depends on $x$ and $\mathcal{F}$.

We define $B_{[\gamma]}^\star:=\{[\tau]\in\mathcal{C} \mid \exists \tau':[0,1]\rightarrow X, \tau'(0)=x,\tau'(0,1]\subset B,\tau=\gamma*\tau'\}$ for every $[\gamma]\in \bar{\varphi}^{-1}(x)$. We claim that $\bar{\varphi}^{-1}(B)=\cup_{\gamma\in\bar{\varphi}^{-1}(x)}B^\star_{[\gamma]}$.
\begin{itemize}[wide]
\item[$(\subset)$] Let $\tau\in\bar{\varphi}^{-1}(B)$, then $\tau(1)\in B$, and we have an $\alpha:[0,1]\rightarrow X$ such that $\alpha(0)=x$, $\alpha(1)=\tau(1)$ and $\alpha(0,1]\subset B$. Observe that $\tau *\alpha^{-1}\in\bar{\varphi}^{-1}(x)$ and $\tau=\tau*\alpha^{-1}*\alpha$. Thus $\tau\in B^\star_{[\tau*\alpha^{-1}]}$.

\item[$(\supset)$] $\tau\in\bigcup_{\gamma\in\bar{\varphi}^{-1}(x)}B^\star_{[\gamma]}$. There exist $\gamma\in\bar{\varphi}^{-1}(x)$ such that $\tau\in B_{[\gamma]}^\star$, and $\alpha:[0,1]\rightarrow X$ such that $\tau=\gamma*\alpha$ and $\alpha(0,1]\subset B$. Thus $\tau(1)\in B$, getting $\tau\in\bar{\varphi}^{-1}(B)$.
\end{itemize}

We verify that $\bar{\varphi}^{-1}$ is homeomorphic to $B\times \bar{\varphi}^{-1}(B)$ observing that this union is disjoint: Let $\gamma,\gamma'\in\bar{\varphi}^{-1}(x)$ such that $B^\star_{[\gamma]}\cap B^\star_{[\gamma']}\neq \varnothing$; let $\tau$ be an element in that intersection. Then there exists $\alpha,\beta:[0,1]\rightarrow X$ such that $\tau=\gamma*\alpha=\gamma*\beta$ and $\alpha(0,1],\beta(0,1]\subset B$. We obtain that $\alpha*\beta^{-1}$ is a closed path, then is contractible, then $[\gamma']=[\gamma*\alpha*\beta^{-1}]=[\gamma]$.

Next phase is to show that $B$ is homeomorphic to $B_{[\gamma]}^\star$ via $\bar{\varphi}|B_{[\gamma]}^\star$. Observe that $\bar{\varphi}|B_{[\gamma]}^\star$ is continuous because $\bar{\varphi}$ also is continuous. It remains to show that it has a continuous inverse.

Let $y\in B$, then there exists $\alpha:[0,1]\rightarrow X$ such that $\alpha(0)=x$, $\alpha(1)=y$ and $\alpha(0,1]\subset B$. Thus $\gamma*\alpha\in B_{[\gamma]}^\star$ and $\gamma*\alpha(1)=y$. Thus, $\varphi|B_{[\gamma]}^\star$ is surjective.

Let $\tau,\tau'\in B_{[\gamma]}^\star$ such that $\tau(1)=\tau'(1)$, then there exits $\alpha,\beta:[0,1]\rightarrow X$ such that $\tau=\gamma*\alpha$, $\tau'=\gamma*\beta$ and $\alpha(0,1],\beta(0,1]\subset B$. $\alpha*\beta^{-1}$ is a closed path in $B$, thus contractible within $X$. Therefore $\tau=\tau'$. We conclude that $\bar{\varphi}|B^\star_{[\gamma]}$ is one-to-one.

It only left to prove that $(\bar{\varphi}|B^\star_{[\gamma]})^{-1}$ is continuous. Let $\mathcal{G}\rightarrow y\in B$ such that $\mathcal{G}=\mathcal{G}\cap [y]$. Since $B$ is path connected and $B\in\mathcal{B}_x$, there exists $\alpha:[0,1]\rightarrow X$ such that $\alpha(0)=x$, $\alpha(1)=y$ and $\alpha(0,1]\subset B$.

Since $\mathcal{G}\rightarrow y$ in $B$, then $[\mathcal{G}]_X\rightarrow y$ in $X$ and there exists $\mathcal{G}'\rightarrow y$ in $X$ such that $\mathcal{G}'\subset [\mathcal{G}]_X$ and there exists a local covering system at $y$ $\mathfrak{D}$ in $X$ such that every element in $\mathfrak{D}$ is connected and $[\mathfrak{D}]_X=\mathcal{G}'$.

For every $D\in\mathfrak{D}$, we define $D^\star_{[\gamma\star\alpha]}\coloneqq \{[\tau]\in \mathcal{C}\mid \exists \tau':[0,1]\rightarrow X,\tau'(0)=y,\tau'(0,1]\subset D,\tau=\gamma*\alpha*\tau'\}$. Observe that $[\gamma * \alpha] \in D^\star_{[\gamma*\alpha]}\cap B^\star_{[\gamma]}$. Hence $[D^\star_{\gamma* \alpha}\mid D\in\mathfrak{D}]_\mathcal{C}\subset [D^\star_{\gamma* \alpha}\cap B^{\star}_{[\gamma]}\mid D\in\mathfrak{D}]_\mathcal{C}\rightarrow [\gamma* \alpha]$ in $\mathcal{C}$. Thus $[D^\star_{\gamma* \alpha}\cap B^{\star}_{[\gamma]}\mid D\in\mathfrak{D}]_{B^\star_{[\gamma]}}\rightarrow [\gamma* \alpha]$ in $B^\star_{[\gamma]}$.

We claim that $[D^\star_{\gamma* \alpha}\cap B^{\star}_{[\gamma]}\mid D\in\mathfrak{D}]_{B^\star_{[\gamma]}}\subset (\bar{\varphi}|B^\star_{[\gamma]})^{-1}(\mathcal{G})$, and thus $(\bar{\varphi}|B^\star_{[\gamma]})^{-1}(\mathcal{G})\rightarrow [\gamma*\alpha]$. Let $D'\in [D^\star_{\gamma* \alpha}\cap B^{\star}_{[\gamma]}\mid D\in\mathfrak{D}]_{B^\star_{[\gamma]}}$, then there exists $D\in\mathfrak{D}$ such that $D^\star_{[\gamma*\alpha]}\cap B^\star_{[\gamma]}\subset D'$. Then $D\in [\mathcal{G}]_X$ and $D\cap B\in\mathcal{G}$. Since $(\bar{\varphi}|B^\star_{[\gamma]})(D^\star \cap B^\star)=D\cap B$, then $(\bar{\varphi}|B^\star_{[\gamma]})^{-1}(D\cap B)\subset D'$, and $D'\in (\bar{\varphi}|B^\star_{[\gamma]})^{-1}(\mathcal{G})$.
\end{proof}

\begin{lema}
\label{lem:CIsSimplyConnected}
$\mathcal{C}$ is path connected and simply connected.
\end{lema}

\begin{proof}
Let $\varepsilon_0$ be the constant path in $x_0$, and let $\gamma\in P(X,x_0)$. For every $u\in [0,1]$, we define the path $\gamma_u:I\rightarrow X$ such that $\gamma_u(t)=\gamma(ut)$. We also define $\tilde{\gamma}:[0,1]\rightarrow \mathcal{C}$ where $\tilde{\gamma}(t)=[\gamma_t]$. If $\tilde{\gamma}$ is continuous, then it is a path between $\varepsilon_0$ and $[\gamma]$.

Let $t\in [0,1]$ and $\mathcal{T}_t:=[(t-\frac{1}{n},t+\frac{1}{n})\mid n\in\mathbb{N}]_I$. Since $\gamma$ is continuous, then $\mathcal{G}:=\gamma(\mathcal{T}_t)\rightarrow \gamma(t)$. We claim that $\tilde{\gamma}(\mathcal{T}_t)=[\gamma(t-\frac{1}{n},t+\frac{1}{n})^\star_{[\gamma_t]}\mid n\in\mathbb{N}]_\mathcal{C}$. It suffices to prove that there exists $N\in\mathbb{N}$ such that for every $n>N$ satisfies $\tilde{\gamma}(t-\frac{1}{n},t+\frac{1}{n})=\gamma(t-\frac{1}{n},t+\frac{1}{n})^\star_{[\gamma_t]}$.

Let $\tau_\varepsilon:[0,1]\rightarrow X$ such that $\tau_\varepsilon(s)=\gamma(s\varepsilon+t)$ for each $\varepsilon\in [-t,1-t]$. Note that $[\gamma_t*\tau_\varepsilon]=[\gamma_{t+\varepsilon}]$. By hypothesis, $X$ is locally path connected and semi-locally simply connected. Then there exists $U\in\mathcal{G}$ path connected whose closed paths are contractible to one point in $X$; consequently, there exists $N\in\mathbb{N}$ such that $\gamma(t-\frac{1}{N},t+\frac{1}{N})\subset U$ has the same properties. $N$ is our candidacy. Let $n>N$
\begin{itemize}[wide]
	\item[$\subset$] Let $\tau\in\tilde{\gamma}(t-\frac{1}{n},t+\frac{1}{n})$, then $\tau=[\gamma_{t+\varepsilon}]$ with $\frac{-1}{n}<\varepsilon<\frac{1}{n}$. By construction, $\tau=[\gamma*\tau_\varepsilon]$ where $\tau_\varepsilon[0,1]\subset \gamma(t-\frac{1}{n},t+\frac{1}{n})$. Thus $\tau\in \gamma(t-\frac{1}{n},t+\frac{1}{n})^\star_{[\gamma_t]}$.
	\item[$\supset$] Let $\tau\in \gamma(t-\frac{1}{n},t+\frac{1}{n})^\star_{[\gamma_t]}$, then there exists $\alpha:[0,1]\rightarrow X$ such that $\tau=\gamma_t*\alpha$ and $\alpha[0,1]\subset \gamma(t-\frac{1}{n},t+\frac{1}{n})$. Define $\varepsilon:=\alpha(1)$. Since $\gamma(t-\frac{1}{n},t+\frac{1}{n})$ satisfies that every closed path within it is contractible in $X$, then $\alpha*\tau_\varepsilon^{-1}\simeq\varepsilon_0$. Thus $\tau=[\gamma_t*\alpha*\tau_\varepsilon^{-1}*\tau_\varepsilon]=[\gamma_{t+\varepsilon}]$ in $\mathcal{C}$.
\end{itemize}
Thus $\tilde{\gamma}$ is continuous. Therefore $\mathcal{C}$ is path connected.

The next phase is to show that $\mathcal{C}$ is simply connected. Let $\bar{\gamma}$ be a closed path at $[\varepsilon_0]$ in $\mathcal{C}$. Observe that $\gamma\coloneqq\bar{\varphi}\bar{\gamma}$ is a closed path at $x_0$ in $X$. For every $u\in I$, we define $\gamma_u: I\rightarrow X$ such $\gamma_u(t)=\gamma(ut)$.

The heart of the proof is to show that $\bar{\gamma}(t)=[\gamma_t]$. Since $\bar{\gamma}$ is continuous, $[\bar{\gamma}(t-\frac{1}{n},t+\frac{1}{n})]_\mathcal{C}$ convergences to $\bar{\gamma}(t)$. Then there exists $\mathcal{F}\rightarrow \bar{\gamma}(t)(1)$ in $X$ such that
\begin{align*}
\left[\bar{\gamma}\left(t-\frac{1}{n},t+\frac{1}{n}\right)\mid n\in\mathbb{N}\right]_\mathcal{C}=[B^\star_{\bar{\gamma}(t)}]_{B\in\mathcal{B}_{\gamma(t)}\cap\mathcal{F}}.
\end{align*}
Also by continuity, $\mathcal{G}\coloneqq[\gamma(t-\frac{1}{n},t+\frac{1}{n})\mid n\in\mathbb{N}]_X\rightarrow \gamma(t)$. Since $X$ is semi-locally simply connected and locally path-connected, there exists $B_t\in\mathcal{B}_{\gamma(t)}\cap\mathcal{F}\cap\mathcal{G}$ connected such that the closed path in $B_t$ are contractible within $X$.

By the equality between the filters above, there exists $N_t\in\mathbb{N}$ such that $\bar{\gamma}(t-\frac{1}{N_t},t+\frac{1}{N_t})\subset (B_t)^\star_{\tilde{\gamma}(t)}$. Thus, for every $\varepsilon\in (\frac{-1}{N_t},\frac{1}{N_t})$ there exists $\alpha_\varepsilon^t$ path in $X$ such that $\alpha_\varepsilon^t(0,1]\subset B_t$ and $\bar{\gamma}(t+\varepsilon)=\bar{\gamma}(t)*[\alpha_\varepsilon^t]$.

On the other hand, $\gamma(t-\frac{1}{N_t},t+\frac{1}{N_t})=\bar{\varphi}\bar{\gamma}(t-\frac{1}{N_t},t+\frac{1}{N_t})\subset B_t$. Observe that $\cup_{t\in [0,1]}(t-\frac{1}{N_t},t+\frac{1}{N_t})=[0,1]$, then there exist $0=t_0<t_1<t_2<\ldots<t_k=1$ such that $\cup_{i=0}^k (t_i-\frac{1}{N_{t_i}},t_i+\frac{1}{N_{t_i}})=[0,1]$ by (topological) compactness. Lastly we observe that:
\begin{itemize}
	\item For $t_0$, $\gamma_\varepsilon[0,1]\subset B_0$ for $\varepsilon<\frac{1}{N_0}$, then $\bar{\gamma}(\varepsilon)*[\gamma_\varepsilon^{-1}]=[\varepsilon_0]$. Thus $\bar{\gamma}(\varepsilon)=\bar{\gamma}(\varepsilon)*[\gamma_\varepsilon^{-1}*\gamma_\varepsilon]=[\gamma_\varepsilon]$.
	\item Without loss of generality, suppose that $[0,1+\frac{1}{N_0})\cap (1-\frac{1}{N_{t_1}},1+\frac{1}{N_{t_1}})\neq\emptyset$. Let $y$ be an element in that intersection. We obtain that
	\begin{align*}
	[\gamma_{t_1}*\tau_{y-t_1}^{t_1}]=[\gamma_y]=\bar{\gamma}(y)=\bar{\gamma}(t_1)*[\alpha_{y-t_1}^{t_1}]
	\end{align*}		
	where $\tau_{y-t_1}^{t_1}$ is path in $X$ such that $\tau_{y-t_1}^{t_1}(s)=\gamma(s(y-t_1)+t_1)$. Thus $[\gamma_{t_1}]=\bar{\gamma}(t_1)$ and satisfies the same for each element in $(t_1-\frac{1}{N_{t_1}},t_1+\frac{1}{N_{t_1}})$.
	\item Repeating that process along $[0,1]$, we obtain that $\bar{\gamma}(t)=\gamma_t$ for every $t\in [0,1]$.
\end{itemize}
Since $\bar{\gamma}(t)=[\gamma_t]$, we observe that $\varepsilon_0=\bar{\gamma}(0)=\bar{\gamma}(1)=[\gamma_1]=[\gamma]\in\mathcal{C}$. We also showed that every closed path in $\mathcal{C}$ is characterized by his image by $\varphi$, thus $\pi(\mathcal{C},\varepsilon_{0})\rightarrow \pi(X,x_0)$ is one-to-one. Therefore, $\bar{\gamma}\simeq 0$.
\end{proof}

\begin{teorema}
$\mathcal{C}$ is the universal cover of $X$.
\end{teorema}

\begin{proof}
\autoref{lem:VarphiIsLocallyTrivialMap} implies that $\varphi:\mathcal{C}\rightarrow X$ is a cover map. \autoref{lem:CIsSimplyConnected} tells us $\mathcal{C}$ is a universal cover of $X$.
\end{proof}

\section{pseudotopological Modification}
\label{cap.ChoquetModification}

Pseudotopological spaces are important because they include graphs, digraphs and scaled metric space, that is, include pretopological spaces. Particularly, pseudotopological spaces are the cartesian closed hull of pretopological spaces (\cite{Beattie_Butzmann_2002}, page 58), what means that they are the smallest cartesian closed category which contains pretopological spaces.

This last section describes the changes we may do to obtain the universal covering space. Despite the fact that most of the definitions are categorical and hence they work for every categorical construct, we use two particular definitions which depends on limit spaces: quotients (\ref{def:Quotients}) and continuous limit structure (\ref{def:ContinuousLimitStructure}).

To develop our objective, we describe an equivalence of quotient in pseudotopological spaces and enunciate results which tell us when a continuous limit structure is a pseudotopological Space. Then we apply that property to prove the last theorems in \autoref{sec: Universal covering space} for our space of interest.

\subsection{Basic notions on pseudotopological Spaces}
\label{sec.BasicNotionsChoquet}

\begin{definicion}
A \emph{pseudotopological Space} (or \emph{Choquet space}) if $\mathcal{F}\rightarrow x$ in $X$ whenever every ultrafilter $\mathcal{G}$ finer than $\mathcal{F}$ converges to $x$ in $X$.
\end{definicion}

There exists a left functor from the pseudotopological spaces (PsTop) to the limit spaces (Lim) called the pseudotopological modifical $\xi:PsTop\rightarrow Lim$ (its right functor is the inclusion $\iota_{PsTop}:Lim\rightarrow PsTop$). A direct consequence of these facts is:

\begin{corolario}[\cite{Beattie_Butzmann_2002}, 1.3.27]
Let $(X_i)_{i\in I}$ be a family of pseudotopological limit spaces and let $X$ be a limit space which carries the initial limit structure with respect to $(f_i:X\rightarrow X_i)_{i\in I}$. Then $X$ is pseudotopological.

In particular subspaces and product of pseudotopological limit spaces are themselves pseudotopological spaces.
\end{corolario}

Since $I$ is a pseudotopological space, homotopy is well-defined in pseudotopological spaces. Also for the continuous limit structure, we have a result which help us in pseudotopological.

\begin{teorema}[\cite{Beattie_Butzmann_2002}, 1.5.5]
Let $X$ and $Y$ be non-empty convergence spaces. Then $\mathcal{C}(X,Y)$ with the continuous limit structure is a pseudotopological space if and only if $Y$ is a pseudotopological space
\end{teorema}

Hence, $P(X,x_0)\coloneqq \{\gamma \in C(I,X):\gamma(0)=x_0\}$ for a fixded $x_0\in X$ as subspace of $C(I,X)$ is a pseudotopological space if $X$ is a pseudotopological space. Thus we only have to modify the category where we made the quotient in the universal cover, observing that in general we do not have that the quotient of any pseudotopological spaces is so (see \autoref{exam:QuotientNotPsTop}). To do that, we remember the definition of quotient through final structures for any topological construct and observe an equivalence in pseudotopological spaces.

\begin{definicion}
Let $\mathcal{D}$ be a topological construct, $X$ be an object in $\mathcal{D}$, $Y$ be a set and $f:X\rightarrow Y$ be a surjective function in $\mathcal{D}$. Then $f$ is a \emph{quotient map} if and only if $Y$ has the final structure with respect to $f$, that is, for any object $Z$ in $\mathcal{D}$ a map $g:Y\rightarrow Z$ is a morphism in $\mathcal{D}$ if and only if the composite map $g\circ f:X\rightarrow Z$ is a morphism in $\mathcal{D}$.
\end{definicion}

\begin{teorema}
Let $X$ be a pseudotopological Space and $Y$ be a set. $f:X\rightarrow Y$ is a quotient map if and only if $\mathcal{F}\rightarrow y$ if for every ultrafilter $\mathcal{G}\supset\mathcal{F}$ there exists $x\in f^{-1}(y)$ and $\mathcal{H}\rightarrow x$ such that $f(\mathcal{H})\subset\mathcal{G}$.
\end{teorema}

\begin{proof}
We start observing that $f:X\rightarrow Y$ is continuous with that structure on $Y$ Consider $\mathcal{F}\rightarrow x$ in $X$, then tautologically $q(\mathcal{F})\subset \mathcal{G}$ for every ultrafilter of $q(\mathcal{F})$ and then $q(\mathcal{F})\rightarrow q(x)$.

Let $g:Y\rightarrow Z$ be a function. Observe that if $g$ is continuous, then $gf$ is continuous by composition of continuous maps. On the other hand, suppose that $gf$ is continuous and $\mathcal{F}\rightarrow y$ in $Y$. Then, for every ultrafilter $\mathbb{G}$ of $\mathcal{F}$ there exists a filter $\overline{\mathcal{G}}\rightarrow x$ such that $f(x)=y$ and $f(\overline{\mathcal{G}})\subset \mathcal{G}$. We define $\mathcal{H}$ as the intersection of all of the ultrafilter of every $\overline{\mathcal{G}}$; $\mathcal{H}\rightarrow x$ because $X$ is a pseudotopological space. We now have that $gf(\mathcal{H})\rightarrow gf(x)=g(y)$.

We claim that $f(\mathcal{H})\subset \mathcal{F}$. Suppose that $A\in f(\mathcal{H})$, then there exists $B\in\mathcal{H}$ such that $f(B)\subset A$, which implies that $B\in\overline{\mathcal{G}}$ for every ultrafilter $\mathcal{G}$ of $\mathcal{F}$ and satisfies that $f(B)\subset A$. Then $f(B)\in \mathcal{G}$ for every ultrafilter $\mathcal{G}$ of $\mathcal{F}$ and satisfies that $f(B)\subset A$. Concluding that $f(B)$ and $A$ are in $\mathcal{F}$. Concluding that $g(\mathcal{F})\rightarrow g(y)$.
\end{proof}

\begin{ejemplo}
\label{exam:QuotientNotPsTop}
The limit quotient of a pseudotopological spaces is not necessary pseudotopological. Consider $I$ with its topological metric structure and define for $x\in I$ and all of the ultrafilters $\mathcal{F}\rightarrow x$ the pseudotopological (actually pretopological) structure on $I$, which we denote as $I_{x,\mathcal{F}}$ such that the only convergent filter in $y\neq x$ is $[y]$ and the only convergent filters in $x$ are $[x]$, $\mathcal{F}$ and $[x]\cap\mathcal{F}$.

Even when the disjoint union is a final structure, the disjoint union of pseudotopological spaces is pseudotopological. And we define the function (as set)
\begin{align*}
\iota:\bigsqcup\limits_{x,\mathcal{F}} I_{x,\mathcal{F}} \rightarrow I \text{ such that }\iota(t,(x,\mathcal{F})=t.
\end{align*}
Observe that $\mathcal{F}\rightarrow x$ in the limit quotient of $\iota$ if and only if $\mathcal{F}$ is finer than the intersection of a finite amount of ultrafilters, and thus it is cleat that topological neighborhood filter of $x$ does not converge to $x$.

This construction in any pseudotopological space $X$ is special. It is called the \emph{ultrafilter modification} of $X$ and provides a right adjoint of the Pseudotopological (\cite{Beattie_Butzmann_2002}, page 58) modification and also provides an example of non-topological compact space (\cite{Beattie_Butzmann_2002}, example 1.4.5).
\end{ejemplo}

\subsection{Covering Spaces on pseudotopological Spaces}
\label{sec.CoveringSpacesChoquetSpaces}

\begin{proposicion}
\label{prop:ContPhiChoquet}
Let $X$ be locally path-connected. $\varphi: P(X,x_0)\rightarrow X$ such that $\gamma\mapsto \gamma(1)$ is surjective and continuous.
\end{proposicion}

\begin{proof}
Since $X$ is locally path-connected, there exists $\gamma_x\in P(X,x_0)$ such that $\gamma_x(1)$ for every $x\in X$, hence $\varphi$ is surjective. Let $\mathcal{H}\rightarrow f$  in $P(X,x_0)$. By continuous limit structure definition, $\omega_{I,X}(\mathcal{H}\times [1])\rightarrow f(1)$. Then $\varphi(\mathcal{H})\rightarrow f(1)$. Thus $\varphi$ is continuous.
\end{proof}

For \ref{prop:AlmostQuotient} and \ref{prop:ContPhiChoquet}, $\varphi$ is a quotient map. Hence we have the lifting theorem \ref{theo:LiftingTheorem} for pseudotopological spaces and if we define $\mathcal{C}\coloneqq \{[\gamma]\mid \gamma\in P(X,x_0)\}$, where $\gamma'\in[\gamma]$ if and only if $\gamma'(1)=\gamma(1)$ and $\gamma\simeq\gamma'$, with quotient pseudotopological structure given by $q:P(X,x_0)\rightarrow \mathcal{C}$ such that $\gamma\mapsto [\gamma]$ and we define $\bar{\varphi}:\mathcal{C}\rightarrow X$ such that $\gamma\mapsto \gamma(1)$, then $\bar{\varphi}$ is continuous and quotient.

To prove the properties of $\bar{\varphi}$ and $\mathcal{C}$ only depends on local connectedness, semi local simply connectednees and $\bar{\varphi}$ being a quotient. Thus $\varphi$ is covering map and $\mathcal{C}$ is simply connected and path-connected. Hence $\bar{\varphi}$ is a universal cover map.

We conclude this section observing that the limit construction of $\mathcal{C}$ for a pseudotopological space $X$ is actually a pseudotopological space. This remark essentially comes from the following result, which tells us that the convergent filters for a path $\gamma$ are very similar in terms of their homotopy.

\begin{lema}
\label{lem:HomFilterMod}
Let $X$ be a connected, locally path connected, semi-locally simply connected pseudotopological space, $\gamma\in P(X,x_0)$ and $\mathcal{F}\rightarrow \gamma$. Then there exists a filter $\mathcal{F}'\rightarrow \gamma$ such that there exist elements $B\in \mathcal{F}$ and $B'\in \mathcal{F}'$ such that
\begin{enumerate}
\item For every $\tau\in B$, there exists $\tau'\in B'$ such that $[\tau] = [\tau']$.
\item Every element in $B'$ is of the form $\gamma*\alpha$.
\end{enumerate}
\end{lema}

\begin{proof}
Since $\mathcal{F}\rightarrow \gamma$, by definition we obtain that $\omega_{I,X}(\mathcal{U}_t\times \mathcal{F})\rightarrow \gamma(t)$ in $X$ for every $t\in I$. We know that $X$ locally path-connected and semi-locally path connected, then there exists $A_t\in \omega_{I,X}(\mathcal{U}_t\times \mathcal{F})$ such that it is path-connected and every loop inside $A_t$ is contractible in $X$. We also obtain the existence of $\varepsilon_t> 0$ and $C_t\in\mathcal{F}$ such that $\omega_{I,X}((t-\varepsilon_t,t+\varepsilon_t)\times C_t) \subset A_t$.

Since $I$ is compact, and we can always assume that $\gamma\in C_t$, we can find a finite amount of real numbers $0=t_0<t_1<\ldots<t_n=1$ such that
\begin{align*}
\bigcup_{\{t_0,t_1,\ldots,t_n\}}C_t \supset \gamma(I)
\end{align*}
We call $C$ the intersection of $C_{t_i}$, observing that it is still in $\mathcal{F}$. Consider $\tau\in C$, then there exists $\alpha:I\rightarrow A_1$ such that $\alpha(0)=\gamma(1)$ and $\alpha(1)=\tau(1)$ for $A_1$ being path-connected. Since every loop in $A_{t_i}$ is contractible in $X$, we have that $\tau\simeq \gamma*\alpha$.

We make this construction for every $\tau\in C$, then we obtain this for every element for every set in the collection $C\cap\mathcal{F}$. For every $D\in C\cap\mathcal{F}$, then we define $\overline{D}$ as the construction above and then we obtain the filter $\mathcal{F}'\coloneqq [\overline{D}\mid D\in C\cap\mathcal{F}]$. Observe that $C$ and $\overline{C}$ satisfies the conditions by construction. 
\end{proof}

\begin{teorema}
\label{Teo:HomFilterMod}
Let $X$ be a connected, locally path connected, semi-locally simply connected pseudotopological space and $\gamma\in P(X,x_0)$. Then $\mathcal{F}'\rightarrow [\gamma]$ if and only if there exist $\mathcal{F}\rightarrow \gamma$ such that $q\mathcal{F}\subset \mathcal{F}'$.
\end{teorema}

\begin{proof}\ \smallskip

$(\Leftarrow)$ This follows directly from $q$ continuous.

$(\Rightarrow)$ Let $\mathcal{F}'\rightarrow [\gamma]$, then there exist a finite amount of $\mathcal{F}_i\rightarrow \gamma_i$ such that $[\gamma_i]=[\gamma]$ and $\mathcal{F}'\supset q(\mathcal{F}_1)\cap\ldots\cap q(\mathcal{F}_n)$. For every $\mathcal{F}_i\rightarrow\gamma_i$ there exists $\mathcal{G}_i\rightarrow\gamma_i$ which elements only depends on $\gamma_i$, and $B_i\in\mathcal{F}_i$ and $B' _i\in\mathcal{G}_i$ as \autoref{lem:HomFilterMod}. Then we might define $\overline{\mathcal{G}}_i$, replacing $\gamma_i$ by $\gamma$, which converges to $\gamma$. Thus $\overline{\mathcal{G}}_1\cap\ldots\cap \overline{\mathcal{G}}_n\rightarrow \gamma$ and $q(\overline{\mathcal{G}}_1\cap\ldots\cap \overline{\mathcal{G}}_n) = q(\mathcal{F}_1)\cap\ldots\cap q(\mathcal{F}_n)\subset \mathcal{F}'$.
\end{proof}

\begin{corolario}
\label{coro:ChoquetCoverSpace}
Let $X$ be a connected, locally path connected, semi-locally simply connected pseudotopological space. If $P(X,x_0)$ is pseudotopological, then $\mathcal{C}$ is pseudotopological. Concluding that the pseudotopological modification of $\mathcal{C}$ is equal to the universal covering space through the pseudotopological quotient.
\end{corolario}

\begin{proof}
Let $\mathcal{F}$ such that all of its ultrafilters converge to $[\gamma]$. Then for every ultrafilter $\mathcal{G}$ of $\mathcal{F}$, there exists $\overline{\mathcal{G}}\rightarrow\gamma$ such that $q\overline{\mathcal{G}}\subset\mathcal{G}$. Define $\mathcal{H}$ as the intersection of every all of the ultrafilters of all of $\overline{\mathcal{G}}$. Given that $P(X,x_0)$ is a pseudotopological space, then $\mathcal{H}\rightarrow \gamma$. Therefore, $q\mathcal{H}\rightarrow [\gamma]$ and $\mathcal{F}\supset q\mathcal{H}$.
\end{proof}

In the same sense of pseudotopological spaces, there exists pretopological (PsTop) and topological (Top) modifications from limit spaces to these categories: $pi$ and $\tau$, respectively. Both functors are also left adjoints (and the inclusions are left adjoints of them). These facts do not tell us anything about the space of continuous functions, but assure that the limit initial structure with respect a (pre)topological is (pre)topological so.

\begin{corolario}
\label{coro:PreTopCoverSpace}
Let $X$ be a connected, locally path connected, semi-locally simply connected pretopological space. If $C(I,X)$ is pretopological, then $\mathcal{C}$ is pretopological.
\end{corolario}

\begin{proof}
For every filter $\mathcal{F}'\rightarrow [\gamma]$ in $\mathcal{C}$ there exists a filter $\mathcal{F}\rightarrow \gamma$ such that $q\mathcal{F}\subset \mathcal{F}'$. Since $P(X,x_0)$ is pretopological, then the neighborhood filter
\begin{align*}
\mathcal{U}_\gamma\coloneqq\bigcap\limits_{\mathcal{F}\rightarrow \gamma}\mathcal{F}\rightarrow \gamma
\end{align*}
Thus we have that
\begin{align*}
q(\mathcal{U}_\gamma) = q\left( \bigcap\limits_{\mathcal{F} \rightarrow \gamma}\mathcal{F} \right)\subset \bigcap\limits_{\mathcal{F} \rightarrow \gamma} q\mathcal{F} = \bigcap\limits_{\mathcal{F}' \rightarrow [\gamma]} \mathcal{F}'=\mathcal{U}_{[\gamma]} \rightarrow [\gamma]
\end{align*}
with $\mathcal{U}_{[\gamma]}$ the neighborhood filter of $[\gamma]$.
\end{proof}

\begin{ejemplo}
$C(I,X)$ is not necessarily pretopological if $X$ is pretopological. To observe this, consider $\mathbb{R}$ with its typical metric space and its closure space given by 
a privileged scale $r>0$, i.e., for every $x\in\mathbb{R}$ we have the neighborhood filter 
\begin{align*}
\mathcal{U}_x\coloneqq [(x-r-\varepsilon,x+r+\varepsilon)\mid \varepsilon>0].
\end{align*}
It is routine to observe that the path
\begin{align*}
f(x) \coloneqq \left\lbrace\begin{array}{ll}
0 & \text{ if } x\in[0,1/2)\\
r & \text{ if } x\in[1/2,1]
\end{array}\right.
\end{align*}
We can also observe that the filters $[f_{x_0}]$ with $f_{x_0}(x)=f(x)$ when $x\neq x_0$ and $f_{x_0}(x)=f(x_0)+r$ converges to $f$ in the continuous limit structure for every $x_0\in I$. However the filter
\begin{align*}
\bigcap_{x_0\in I}[f_{x_0}] = [f_{x_0}\mid x_0\in I] \supset \mathcal{U}_f
\end{align*}
does not converges to $f$.
\end{ejemplo}

The following proposition tell us that for every topological space $X$, $C(I,X)$ is always topological given the good properties of $I$. This fact, in addition that every limit subspace of topological spaces is a topological space, gives as a corollary (\ref{coro:TopCoverSpace}) that this constructions is actually an extension of the cover spaces in Top and not just an imitation  of the techniques.

\begin{proposicion}
\label{prop:C(I,X)topifXtop}
Let $X$ be a topological space then the continuous limit structure in $C(I,X)$ is continuous, in particular it has the compact-open topological structure.
\end{proposicion}

\begin{corolario}
\label{coro:TopCoverSpace}
Let $X$ be a connected, locally path connected, semi-locally simply connected topological space. Then $\mathcal{C}$ is pretopological. Concluding that the pretopological modification of $\mathcal{C}$ is equal to the universal covering space through the pseudotopological quotient.
\end{corolario}

\begin{proof}
Given \autoref{prop:C(I,X)topifXtop}, we know that $P(X_0,x_0)$ is pretopological and then $\mathcal{C}$ is so. Thus it is enough to show that for every $U$ in the neighborhood filter $\mathcal{U}_{[\gamma]}$ of $[\gamma]$, there exists $V\in\mathcal{U}_{[\gamma]}$ such that for every $[\tau]\in V$ the set $U$ is in the neighborhood filter of $[\tau]$. By construction, observe thatwe can find a $U'$ in $\mathcal{U}_\gamma$ such that $q(U')=U$ and, since $X$ is topological, there exists $V'$ in $\mathcal{U}_\gamma$ such that every $\tau\in V'$ satisfies that $U'\in \mathcal{U}_{\tau}$. Therefore $q(V')\in\mathcal{U}_{[\gamma]}$ and every $[\tau]\in q(V')$ satisfies that $U\in \mathcal{U}_{[\tau]}$.
\end{proof}

\section*{Acknowledgments}

The author would like to thank Omar Antolín y Antonio Rieser for the discussions held regarding this article during the stay in Florida, as well as the suggestions received on subsequent occasions from Antonio Rieser, which is my advisor during my PhD studies. Additionally, the author also thanks the opportunity to participate in that event in Cuernavaca in 2022 supported by the Grant No. DMS-1928930.

\bibliographystyle{amsalpha}
\bibliography{all}


\end{document}